\documentclass[12pt]{amsart}

\usepackage[latin1]{inputenc}

\usepackage{multirow} 
\usepackage{subfigure}
\usepackage{pdflscape}
\usepackage{amssymb,amsmath,epsfig}
\usepackage{float}

\usepackage{hyperref}

\newtheorem{Theorem}{Theorem}[section]
\newtheorem{Corollary}[Theorem]{Corollary}
\newtheorem{Proposition}[Theorem]{Proposition}

\numberwithin{equation}{section}

\newcommand{\N}{{\mathbb N}}

\def\qed{\hfill$\Box$\\ \medskip}

\def\ds{\displaystyle}    
\def\ov{\overline}

\def\a{{\bf a}}
\def\x{{\bf x}}
\def\y{{\bf y}}

\def\1{{\frac{1}{2}}}
\def\3{{\frac{3}{2}}}
\def\5{{\frac{5}{2}}}

\def\ch{{\rm ch\,}}

\def\mod{\,{\rm mod}\,}
\def\sgn{\,{\rm sgn}\,}

\def\ovx{\ov{x}}
\def\ovy{\ov{y}}

\def\cal{\mathcal}

\title{Generating functions for some series of characters of classical Lie groups}

\author{Ronald C. King}

\thanks{
School of Mathematical Sciences, University of Southampton, 
Southampton SO17 1BJ, England ({\tt r.c.king@soton.ac.uk})}

\begin{document}

\begin{abstract}
There exist a number of well known multiplicative generating functions for series of Schur functions. Amongst these
are some related to the dual Cauchy identity whose expansion coefficients are rather simple, and in some cases
periodic in parameters specifying the Schur functions. More recently similar identities have been found involving expansions in
terms of characters of the symplectic group. Here these results are extended and generalised to all classical Lie groups. 
This is done through the derivation of explicit recurrence relations for the expansion coefficients
based on the action of the Weyl groups of both the symplectic and orthogonal groups. Copious results are tabulated in the form of 
explicit values of the expansion coefficients as functions of highest weight parameters. An alternative approach is then based on 
dual pairs of symplectic and/or orthogonal groups. A byproduct of this approach is that expansions in terms of spin orthogonal 
group characters can always be recovered from non-spin cases. 
\end{abstract}

\maketitle

\section{Introduction}\label{sec-introduction}
This article has been prompted by the work of 
Lee and Oh in~\cite{LO1} and~\cite{LO2} in which they studied auto-correlation functions and their evaluation
in terms of characters of irreducible representations of the symplectic and unitary groups, respectively. 
In doing so they arrived at a variety of identities in the form of rather neat generating functions for particular 
series of these characters.

For example, for any $n\in\N$ and indeterminates $\x=(x_1,x_2,\ldots,x_n)$, with inverses $\ov\x=(x_1^{-1},x_2^{-1},\ldots,x_n^{-1})$, 
they established in~\cite{LO1} an identity that can be expressed in the form
\begin{equation}\label{eqn-spn-LO1ex}
   \prod_{i=1}^n (1+x_i^2+x_i^{-2}) = \sum_{r=0}^n \sum_{k=0}^{\lfloor \frac{n-p}{2} \rfloor}\ \psi_{k,r}\ \ch^{Sp(2n)}_{(2^p,1^{2k})}(\x,\ov\x)\,,
\end{equation}
where $r=n-p-2k$, and $\psi_{k,r}$ depends only on $k \mod3$ and $r \mod6$ with $\psi_{k,r}$ given for $k\in\{0,1,2\}$ and $r\in\{0,1,2,3,,4,5\}$ by
\begin{footnotesize}
\begin{equation}\label{eqn-psi3}
\begin{array}{|c||c|c|c|c|c|c|}
                      \hline
                      k\backslash r&0&1&2&3&4&5\cr
                      \hline\hline
                      0&1&0&0&-1&0&0\cr
                       \hline
											1&-1&1&0&1&-1&0\cr
                       \hline
                      2&0&-1&0&0&1&0\cr
											\hline
\end{array}
\end{equation}
\end{footnotesize}

Similarly in~\cite{LO2} they derived a unitary group identity which applies equally well to the general linear group in the form
\begin{equation}\label{eqn-gln-LO2ex}
    \prod_{i=1}^n (1+x_i)(1+x_i^2)=\sum_{p=0}^n \sum_{q=0}^{n-p} \sum_{r=0}^{n-p-q}\, \tau_{q,r}\, \ch^{GL(n)}_{(3^p,2^q,1^r)}(\x)\,,
\end{equation}
where $\tau_{q,r}$ depends only on $q$ and $r \mod 4$ and $\tau_{q,r}$ is given for $q,r\in\{0,1,2,3\}$ by
\begin{footnotesize}
\begin{equation}\label{eqn-tau}
\begin{array}{|c||c|c|c|c|}
\hline
q\backslash r&0&1&2&3\cr
\hline\hline
0&1&1&0&0\cr
\hline
1&1&0&-1&0\cr
\hline
2&0&-1&-1&0\cr
\hline
3&0&0&0&0\cr
\hline
\end{array}
\end{equation}
\end{footnotesize}

If we introduce extra parameters $\y=(y_1,y_2,\ldots,y_m)$ with inverses $\ov\y=(y_1^{-1},y_2^{-1},\ldots,y_m^{-1})$,
these identities can be viewed as particular specialisations of a version of the dual Cauchy identity~\cite{Mac,BG}
\begin{equation}\label{eqn-gln-glm}
     \prod_{i=1}^n\prod_{j=1}^m (x_i+y_j) = \sum_{\lambda\in{m^n}} \ch_{\lambda}^{GL(n)}(\x)\ \ch_{\tilde\lambda}^{GL(m)}(\y) 
\end{equation}
and its symplectic analogue~\cite{BG}
\begin{equation}\label{eqn-spn-spm}
     \prod_{i=1}^n\prod_{j=1}^m (x_i+x_i^{-1}+y_j+y_j^{-1}) = \sum_{\lambda\in{m^n}} \ch_{\lambda}^{Sp(2n)}(\x,\ov\x)\ \ch_{\tilde\lambda}^{Sp(2m)}(\y,\ov\y)\,,
\end{equation}
where in both cases the summation on the right is over partitions $\lambda=(\lambda_1,\lambda_2,\ldots,\lambda_n)$ with no more 
than $n$ nonvanishing parts and largest part $\lambda_1\leq m$, and $\tilde\lambda=(n-\lambda'_m,\ldots,n-\lambda'_2,n-\lambda'_1)$
where $\lambda'=(\lambda'_1,\lambda'_2,\ldots,\lambda'_m)$ is the partition conjugate to $\lambda$.

If one sets $m=3$ and $\y=(1,i,-i)$ in (\ref{eqn-gln-glm}) then the left hand side reduces to that of (\ref{eqn-gln-LO2ex}).
In such a case $\lambda$ is necessarily of the form $(3^p,2^q,1^r)$ and $\tilde\lambda=(s+q+r,s+r,s)$, with $s=n-p-q-r\geq0$,
so that on the right hand side of (\ref{eqn-gln-LO2ex}) one finds that $\tau_{q,r}=\ch_{(q+r,r)}^{GL(3)}(1,i,-i)$, where 
use has been made of the fact that $\ch_{(s+q+r,s+r,s)}^{GL(3)}(1,i,-i)=\ch_{(q+r,r)}^{GL(3)}(1,i,-i)$ for all $s$. 
For given $(q,r)$ the coefficient $\tau_{q,r}$ may then be evaluated from the Weyl character formula for $GL(3)$~\cite{Mac,BG,FH}.
This yields the values of $\tau_{q,r}$ as given above with the periodic behaviour a consequence of the components of $\y$ 
being powers of $i$. 
 
Similarly, if one sets $m=2$ and $\y=(\omega,-\omega)$ with $\omega=e^{i\pi/3}$ in (\ref{eqn-spn-spm}) then the left hand side 
coincides with that of (\ref{eqn-spn-LO1ex}). This time $\lambda$ is necessarily of the form $(2^p,1^q)$ with $\tilde\lambda=(q+r,r)$
with $r=n-p-q\geq 0$. In order to arrive at the right hand side of (\ref{eqn-spn-LO1ex}) it is then necessary to evaluate
$\ch_{(q+r,r)}^{Sp(4)}(\omega,-\omega,\omega^{-1},-\omega^{-1})$ explicitly. Weyl's character formula for $Sp(4)$~\cite{Mac,BG,FH}
implies that these characters are zero unless $q$ is even, and for $q=2k$ yields the values of $\psi_{k,r}$ as given above
with the periodicity this time a consequence of the components of $(\y,\ov\y)$ being powers of $\omega=e^{i\pi/3}$.

In the case of $GL(n)$ a rather different approach to identities of the type (\ref{eqn-gln-LO2ex}) has been offered in work 
on Schur functions that appeared in the mathematical physics literature in the 1980's~\cite{YW,LP,KWY}. In particular the 
identity (\ref{eqn-gln-LO2ex}) is given explicitly in~\cite{KWY}, along with a number of other identities obtained in~\cite{LO2}. 
In this approach, the use of the dual Cauchy identity (\ref{eqn-gln-glm}) is replaced by the following:
\begin{equation}\label{eqn-akappa-lambda}
 \prod_{i=1}^n \left(\sum_{k=0}^m\, a_k x_i^k\,\right) = \sum_{\lambda\in(m^n)} \sum_{\kappa}\,  a(\kappa)\, \epsilon(\kappa,\lambda)\, \ch_\lambda^{GL(n)}(\x)\,,
\end{equation}
where $a(\kappa)=a_{\kappa_1}a_{\kappa_2}\cdots a_{\kappa_n}$ is the coefficient of $\x^\kappa=x_1^{\kappa_1}x_2^{\kappa_2}\cdots x_n^{\kappa_n}$ in the expansion of the product on the left, and the passage from the integer sequence $\kappa$ to the partition $\lambda$ is mediated by a sequence of modifications of the form
\begin{equation}\label{eqn-mu-modification}
     (\mu_1,\ldots,\mu_i,\mu_{i+1}\ldots,\mu_n) \rightarrow (\mu_1,\ldots,\mu_{i+1}-1,\mu_i+1,\ldots,\mu_n)
\end{equation}
for $i=1,2,\ldots,n-1$, with each such modification in the passage from $\kappa$ to $\lambda$ contributing a factor 
of $-1$ to $\epsilon(\kappa,\lambda)$, and with $\epsilon(\kappa,\lambda)=0$ in any case for which a sequence $\mu$ is  encountered 
in which $\mu_i=\mu_{i+1}-1$ for any $i$. 

A consequence of this approach is that it rather readily gives rise to recurrence relations for the coefficients on the right. For example,
setting $m=3$ and $a_0=a_1=a_2=a_3=1$ the left hand side of (\ref{eqn-akappa-lambda}) coincides with that of (\ref{eqn-gln-LO2ex}). The 
general recurrence relation to be found in~\cite{KWY} is such that on the right hand side of (\ref{eqn-gln-LO2ex}) one then has 
coefficients given by:
\begin{equation}\label{eqn-tau-QR}
\begin{array}{l}
\tau_{q,r}=Q_qR_r-Q_{q-1}R_{r-1}~~\mbox{with}~~Q_t=R_t=T_t,\cr
\mbox{where}~~T_0=1, T_1=1, T_2=0~~\mbox{and}~~T_t=T_{t-1}-T_{t-2}+T_{t-3}~~\mbox{for $t\geq3$}\,.\cr
\end{array}
\end{equation}
In this particular example, these recurrence relations led to the explicit multiplicity free expansion of the left hand side of 
(\ref{eqn-gln-LO2ex}) first given in equation (5.15) of~\cite{KWY}, and subsequently rederived in~\cite{LO2} and expressed in the form
(\ref{eqn-gln-LO2ex}) and (\ref{eqn-tau}). 

Many other results of the same type are provided in the three papers~\cite{YW,LP,KWY} including the replacement of $(1+x_i)(1+x_i^2)$
by $(1+x_i)(1-x_i^2)$, $(1-x_i)(1+x_i^2)$ and $(1-x_i)(1-x_i^2)$ in~\cite{KWY}, and rather more simply by $(1\pm x_i^p)$ 
for $1\leq p\leq 4$ in~\cite{YW} and for all $p$ in~\cite{LP}. Moreover the generalisation of the recurrence relations (\ref{eqn-tau-QR})
to the case $\a=(1,a,b,c)$ has been given in~{KWY}.

The identity (\ref{eqn-akappa-lambda}) can be seen as arriving directly from the symmetry of left hand side with 
respect to permutations of the components of $\x$ and the properties of symmetric functions. These observations
lead to the identity:
\begin{equation}\label{eqn-amu-lambda}
 \prod_{i=1}^n \left(\sum_{k=0}^m\, a_k x_i^k\,\right) = \sum_{\mu\in(m^n)} a(\mu)\,m_\mu(\x) 
           =\sum_{\mu\in(m^n)}\sum_{\lambda\in(m^n)}\,a(\mu)\,K^{-1}_{\mu\lambda}\,\ch_\lambda^{GL(n)}(\x)\,,
\end{equation}
where $\mu$, like $\lambda$, is summed over all partitions into no more that $n$ parts with largest part no larger than $m$, 
while $K^{-1}$ is the inverse of the Kostka matrix, that is the transition matrix from the monomial symmetric functions $m_\mu(\x)$ to the 
Schur functions $s_\lambda(\x)=\ch_\lambda^{GL(n)}(\x)$~\cite{Mac}. The equivalence of this formula to (\ref{eqn-akappa-lambda})
can be seen by expressing $m_\mu(\x)$ as a sum of all distinct monomials $\x^\kappa$ with $\kappa$ a permutation of the parts of 
the partition $\mu$, and by evaluating $K^{-1}_{\mu\lambda}$ in terms of signed rim hooks or special border strips as spelled out in~\cite{ER,Mac}.
These rim hooks or strips are nothing other than the slinkies used in~\cite{KWY}. They have lengths $\kappa_i$ and here they
are weighted by $a_{\kappa_i}$ for $i=1,2,\ldots,n$, together with sign factors that are automatically encapsulated in the 
modification rules (\ref{eqn-mu-modification}).

In what follows the approach to character identities like that for $GL(n)$ based on (\ref{eqn-akappa-lambda}) is extended to
all the classical Lie groups, $GL(n)$, $Sp(2n)$, $SO(2n+1)$ and $SO(2n)$. The key to doing this is laid down in 
Section~\ref{sec-char-G} in which the Weyl group symmetry of the product on the left of each identity is exploited 
to give a general formula for the sum on the right expressed in terms of the sign factors of the Weyl group elements
whose action maps vectors $\kappa$ in the weight space, $\Lambda^G$, of each group $G$ to vectors $\lambda$ in $\Lambda^G_+$,
the dominant chamber of the weight space. In each case the required action of the $n$ generators of the relevant Weyl group 
is tabulated

As a precursor to dealing with the symplectic and orthogonal groups, the case of $GL(n)$ treated in~\cite{KWY} and~\cite{LO2}
is recast in Section~\ref{sec-gln} in the more general setting provided by Section~\ref{sec-char-G}. The systematic 
use of the so-called dot action of the Weyl group generators is used to derive the recurrence relations 
generalising those of (\ref{eqn-tau-QR}) in the case of the product over $i$ from $1$ to $n$ 
of the sum $a_0+a_1x_i+a_2x_i^2+a_3x_i^3$ for arbitrary parameters $\a=(a_0,a_1,a_2,a_3)$, along with their successive restriction 
to $a_3=0$ and $a_2=a_3=0$. It is pointed out that the dual Cauchy identity provides an alternative method of 
calculating the expansion coefficients appearing in such $GL(n)$ identities, without however leading so easily to recurrence relations
for these coefficients.

The same process is then repeated for $Sp(2n)$, $SO(2n+1)$ and $SO(2n)$ in Sections~\ref{sec-spn}, \ref{sec-oon} and \ref{sec-eon},
respectively, for the product over $i$ from $1$ to $n$ of the sum $a_0+a_1(x_i+x_i^{-1})+a_2(x_i^2+x_i^{-2})$ for arbitrary $\a=(a_0,a_1,a_2)$,
together with their restriction first to $a_2=0$ and then $a_1=0$. 
The same approach is applied to the expansion of the product over $i$ from $1$ to $n$ of 
$a_\1(x_i^\1+x_i^{-\1})+a_3(x_i^\3+x_i^{-\3})+a_\5(x_i^\5+x_i^{-\5})$ in terms of spin characters of $SO(2n+1)$ and $O(2n)$,
first in the case $a_\5=0$ in Section~\ref{sec-spin} and later in the considerably more complicated case $a_\3=0$ in Appendix~\ref{Appendix}.
Dual pair groups~\cite{Bau,Mor} are exploited in Section~\ref{sec-dual-pair} to provide explicit character theoretic formulae for 
all of the expansion coefficient not only in the symplectic case, but also in both ordinary and spin orthogonal group cases.
These formulae are also shown to provide a precise connection between spin and non-spin character expansion coefficients 
that is exemplified in the $m=1$ and $m=2$ cases.
 
Explicit results corresponding to solutions of the recurrence relations are then tabulated for various $m=1$ and $m=2$ cases 
in the two subsections of Section~\ref{sec-spn-oon-eon-values}. An emphasis is placed on those that are periodic or expressible 
somewhat simply in the spirit of those already identified in~\cite{KWY,YW} and by Lee and Oh in~\cite{LO2,LO1}. 
A few concluding remarks are made in Section~\ref{sec-conclusion}.


\section{Characters of classical Lie groups}\label{sec-char-G} 
Let $n\in\N$ be fixed and let ${\cal E}_n$ be an $n$-dimensional Euclidean space with basis 
provided by mutually orthogonal unit vectors $\epsilon_i$ for $i=1,2,\ldots,n$. Any vector
$\kappa\in{\cal E}_n$ given in this basis by $\kappa=\kappa_1\epsilon_1+\kappa_2\epsilon_2+\cdots+\kappa_n\epsilon_n$
is denoted here more simply by $\kappa=(\kappa_1,\kappa_2,\ldots,\kappa_n)$. 
It is also convenient to denote $x_1^{\kappa_1}x_2^{\kappa_2}\cdots x_n^{\kappa_n}$ by $\x^\kappa$  
for any sequence of indeterminates $\x=(x_1,x_2,\ldots,x_n)$. 
 
Let $G$ be any one of the classical Lie groups $GL(n)$, $SO(2n+1)$, $Sp(2n)$ or $SO(2n)$.
All of these groups possess finite-dimensional irreducible representations $V_G^\lambda$ 
specified by their highest weight $\lambda$ for any $\lambda\in\Lambda_G^+$, the subset of dominant weights
of the weight lattice $\Lambda_G$ of $G$ which in each case can be embedded in ${\cal E}_n$.
The character of $V_G^\lambda$ is given by Weyl's formula
\begin{equation}\label{eqn-weyl-char}
   \ch^G_\lambda(\x) = \sum_{w\in W_G} \sgn(w)\,\x^{w(\lambda+\rho_G)} \bigg{/} \sum_{w\in W_G} \sgn(w)\,\x^{w(\rho_G)} \,, 
\end{equation}  
where $W_G$ is the Weyl group of $G$, whose elements $w$ have signature or parity $\sgn(w)=\pm1$, 
while $\rho_G$ is the Weyl vector of $G$, that is half the sum 
of the positive roots of the corresponding Lie algebra, and $\x=(e^{\epsilon_1},e^{\epsilon_2},\ldots,e^{\epsilon_n})$,
where $e^{\epsilon_i}$ for $i=1,2,\ldots,n$ are formal exponentials of the orthonormal basis vectors of the 
$n$-dimensional space ${\cal E}_n$ in which the weights are embedded.

One may conveniently extend the domain of the right hand side of (\ref{eqn-weyl-char}) to define 
\begin{equation}\label{eqn-char-kappa}
    \ch^G_\kappa(\x) = \sum_{w\in W_G} \sgn(w)\,\x^{w(\kappa+\rho_G)} \bigg{/} \sum_{w\in W_G} \sgn(w)\,\x^{w(\rho_G)} \,, 
\end{equation}
for any $\kappa\in\Lambda_G$, and any sequence of indeterminates $\x=(x_1,x_2,\ldots,x_n)$. It should be noted that
for any $w\in W_G$ we have
\begin{equation}\label{eqn-ch-wkappa}
    \ch^G_\kappa(\x) = \sgn(w)\ \ch^G_{w(\kappa+\rho_G)-\rho_G}(\x)\,.  
\end{equation}

Any Weyl group invariant expression linear in terms of the form $\x^\kappa$ with $\kappa\in\Lambda_G$
can itself be expressed as a linear sum of irreducible characters $\ch^G_\lambda(\x)$ with $\lambda\in\Lambda_G^+$
evaluated at $\x$, that is to say we have: 
\begin{Proposition}\label{Pro-WG-invariance}
For all coefficients $a(\kappa)$ such that
\begin{equation}\label{eqn-weyl-invariant}
    \sum_{\kappa\in\Lambda_G} a(\kappa)\,\x^\kappa =\sum_{\kappa\in\Lambda_G} a(\kappa)\,\x^{w(\kappa)}
\end{equation}	
for any $w\in W_G$, we have 
\begin{equation}\label{eqn-kappa-lambda}
    \sum_{\kappa\in\Lambda_G} a(\kappa)\,\x^\kappa  = \sum_{\lambda\in\Lambda^+_G} \sum_{\kappa\in\Lambda_G} \sgn(w_{\kappa,\lambda})\,a(\kappa)\,\ch^G_\lambda(\x)\,.
\end{equation}
where the summation is over those $\kappa$ for which there exists $w_{\kappa,\lambda}\in W_G$ 
such that $w_{\kappa,\lambda}(\kappa+\rho_G)-\rho_G=\lambda$.
\end{Proposition}

\noindent{\bf Proof}:~~
It follows from (\ref{eqn-weyl-invariant}) that
\begin{equation}	
\begin{array}{l}
		\big(\ds \sum_{\kappa\in\Lambda_G} a(\kappa)\x^\kappa\big) \big(\ds  \sum_{w\in W_G} \sgn(w)\x^{w(\rho_G)}\big)
		= \ds\sum_{w\in W_G} \sgn(w)\big(\ds \sum_{\kappa\in\Lambda_G} a(\kappa)\x^{w(\kappa)}\big) \x^{w(\rho_G)} \cr
		= \ds\sum_{w\in W_G} \sgn(w)\big(\ds \sum_{\kappa\in\Lambda_G} a(\kappa)\x^{w(\kappa+\rho_G)}\big) 
		= \ds\sum_{\kappa\in\Lambda_G} a(\kappa)\big(\ds \sum_{w\in W_G} \sgn(w)\x^{w(\kappa+\rho_G)}\big),\cr 
\end{array}
\end{equation}
so that
\begin{equation}\label{eqn-xkappa-chkappa}
		\sum_{\kappa\in\Lambda_G} a(\kappa)\,\x^\kappa = \sum_{\kappa\in\Lambda_G} a(\kappa)\,\ch^G_\kappa(\x)\,.
\end{equation}
However, if there exists $w\in W_G$ such that $w(\kappa+\rho_G)=\kappa+\rho_G$ with $\sgn(w)=-1$ then $\ch^G_\kappa(\x)=0$, while in all other
cases for each $\kappa\in\Lambda_G$ there exists a unique dominant weight $\lambda\in\Lambda^+_G$ and some $w_{\kappa,\lambda}\in W_G$ of signature $\sgn(w_{\kappa,\lambda})$ such that $w_{\kappa,\lambda}(\kappa+\rho_G)=\lambda+\rho_G$, and in such a case $\ch^G_\kappa(\x)=\sgn(w_{\kappa,\lambda})\,\ch^G_\lambda(\x)$.
This suffices to complete the proof of (\ref{eqn-kappa-lambda}).
\qed

The data we require on the Weyl groups $W_G$ of $G$, 
that is their sets of positive roots, $\Delta_G^+$, together with half their sum, $\rho_G$,
are given in Table~\ref{Tab-roots-rho}~\cite{Hum,FH}. 

\begin{table}[ht]
\begin{center}
\begin{tabular}{|l|l|l|}
\hline
$G$&Positive roots~~$\Delta_G^+$&$\rho_G=\frac{1}{2}\sum_{\alpha\in\Delta_G^+}\alpha$\cr
\hline\hline$GL(n)$&$\{\epsilon_i-\epsilon_j|1\leq i<j\leq n\}$&$\sum_{i=1}^n (n-i)\epsilon_i$\cr
\hline
$SO(2n+1)$&$\{\epsilon_i\pm\epsilon_j|1\leq i<j\leq n\}\cup\{\epsilon_i|1\leq i\leq n\}$&$\sum_{i=1}^n (n+\frac{1}{2}-i)\epsilon_i$\cr
\hline
$Sp(2n)$&$\{\epsilon_i\pm\epsilon_j|1\leq i<j\leq n\}\cup\{2\epsilon_i|1\leq i\leq n\}$&$\sum_{i=1}^n (n+1-i)\epsilon_i$\cr
\hline
$SO(2n)$&$\{\epsilon_i\pm\epsilon_j|1\leq i<j\leq n\}$&$\sum_{i=1}^n (n-i)\epsilon_i$\cr
\hline
\hline
\end{tabular}
\end{center}
\medskip
\caption{The positive roots $\alpha\in\Delta_G^+$ of the classical Lie groups $G$ and half their sum, $\rho_G$.} 
\label{Tab-roots-rho}
\end{table}

The corresponding Weyl groups $W_G$ may be identified as the symmetric group $S_n$ in the case of $GL(n)$,
as the semidirect product $(Z_2)^n\rtimes S_n$ in the case of both $SO(2n+1)$ and $Sp(2n)$, and
as $(Z_2)^{n-1}\rtimes S_n$ in the case of $SO(2n)$~\cite{Hum}. In its action on $\kappa$ the group $S_n$ 
is the group of permutations of $(\kappa_1,\kappa_2,\ldots,\kappa_n)$, while $(Z_2)^n\rtimes S_n$ is the group of all permutations 
and sign changes of the components of $\kappa$ and $(Z_2)^{n-1}\rtimes S_n$ is its subgroup involving an even number of sign changes.

The action of $W_G$ is generated by reflections $w_\alpha$ in the hyperplanes perpendicular to the simple roots $\alpha\in\Pi_G$. 
Their action to give $w_\alpha(\kappa+\rho_G)-\rho_G$ for any $\kappa=(\kappa_1,\kappa_2,\ldots,\kappa_n)$ is given
explicitly in Table~\ref{Tab-wkappa} in which any components indicated by $\ldots$ are left unchanged under the reflections. 

\begin{table}[ht]
\begin{center}
\begin{tabular}{|l|l|l|}
\hline
$G$&Simple roots~~$\alpha\in\Pi_G$&$\begin{array}{l} w_\alpha(\kappa+\rho_G)-\rho_G~~\mbox{with}\cr
                       \kappa=(\kappa_1,\ldots,\kappa_i,\kappa_{i+1},\ldots,\kappa_n)\cr
											\end{array}$\cr
\hline\hline$GL(n)$&$\alpha_i=\epsilon_i-\epsilon_{i+1}$ with $1\leq i<n$&$(\ldots,\kappa_{i+1}-1,\kappa_{i}+1,\ldots)$\cr
\hline
$SO(2n+1)$&$\alpha_i=\epsilon_i-\epsilon_{i+1}$ with $1\leq i<n$&$(\ldots,\kappa_{i+1}-1,\kappa_{i}+1,\ldots)$\cr
             &$\alpha_n=\epsilon_{n}$&$(\ldots,-\kappa_{n}-1)$\cr
\hline
$Sp(2n)$&$\alpha_i=\epsilon_i-\epsilon_{i+1}$ with $1\leq i<n$&$(\ldots,\kappa_{i+1}-1,\kappa_{i}+1,\ldots)$\cr
           &$\alpha_n=2\epsilon_{n}$&$(\ldots,-\kappa_{n}-2)$\cr
\hline
$SO(2n)$&$\alpha_i=\epsilon_i-\epsilon_{i+1}$ with $1\leq i<n$&$(\ldots,\kappa_{i+1}-1,\kappa_{i}+1,\ldots)$\cr
           &$\alpha_n=\epsilon_{n-1}+\epsilon_{n}$&$(\ldots,-\kappa_{n}-1,-\kappa_{n-1}-1)$\cr
\hline
\hline
\end{tabular}
\end{center}
\medskip
\caption{The simple roots $\alpha\in\Pi$ of the classical Lie groups $G$ and their action $w_\alpha(\kappa+\rho_G)-\rho_G$.} 
\label{Tab-wkappa}
\end{table}


\section{General linear group character identities}\label{sec-gln}
Let $n\in\N$ be fixed and let $\x=(x_1,x_2,\ldots,x_n)$, and let $\ch_{GL(n)}^\lambda(\x)$ denote the character 
of the irreducible representation of $GL(n)$ of highest weight $\lambda$ evaluated on a group element with eigenvalues $\x$.

\begin{Theorem}[\cite{KWY} see (5.11) for the case $a_0\!=\!1$]\label{The-gln}
For all $\a\!=\!(a_0,a_1,a_2,a_3)$ we have
\begin{equation}\label{eqn-gln-gf3}
 \prod_{i=1}^n (a_0+a_1x_i+a_2x_i^2+a_3x_i^3) =  \sum_{p=0}^n\ \sum_{q=0}^{n-p}
             \sum_{r=0}^{n-p-q}\ a_3^p\,\psi_{q,r}(\a)\,a_0^{n-p-q-r}\ \ch^{GL(n)}_{(3^p,2^q,1^r)}(\x)\,,
\end{equation}
where
\begin{equation}\label{eqn-gln-psi}
     \psi_{q,r}(\a) = 	Q_q\,R_r-a_0a_3\,Q_{q-1}\,R_{r-1}\,,
\end{equation}
with $Q_q=0$ and $R_r=0$ if $q<0$ and $r<0$, respectively, while
\begin{equation}\label{eqn-gln-Q}
Q_0=1~~\mbox{and}~~Q_q=a_2\,Q_{q-1}-a_1a_3\,Q_{q-2}+a_0a_3^2\,Q_{q-3}~~\mbox{for $q\geq1$}\,,
\end{equation} 
and
\begin{equation}\label{eqn-gln-R}
R_0=1~~\mbox{and}~~R_r=a_1\,R_{r-1}-a_0a_2\,R_{r-2}+a_0^2a_3\,R_{r-3}~~\mbox{for $r\geq1$}\,.
\end{equation}
\end{Theorem}

\noindent{\bf Proof}:~~
A proof has been provided in the $n$-independent case $\a=(1,a,b,c)$ in~\cite{KWY}, 
with the expansion coefficients given by $g_{pqr}(abc)$ in (5.11). The full result for all $a_0$
may then be recovered merely by exploiting homogeneity in the total powers of the various $a_i$. However,
the derivation is reproduced here in a manner more obviously appropriate for extension
to the case of characters of the other classical Lie groups, both orthogonal and symplectic.

First we note that
\begin{equation}\label{eqn-gln-gf3-xkappa}
   \prod_{i=1}^n (a_0+a_1x_i+a_2x_i^2+a_3x_i^3)=\sum_{\kappa}  a_3^k\,a_2^\ell\,a_1^m\,a_0^{n-k-\ell-m} \x^\kappa
\end{equation}
where $\kappa_j\in\{0,1,2,3\}$ for $j=1,2,\ldots,n$ and
\begin{equation}
\begin{array}{rcl}
     k&=&\#\{\kappa_j=3|j=1,2,\ldots,n\},\cr
  \ell&=&\#\{\kappa_j=2|j=1,2,\ldots,n\},\cr
		 m&=&\#\{\kappa_j=1|j=1,2,\ldots,n\}.\cr
\end{array}
\end{equation}

The left hand side of (\ref{eqn-gln-gf3-xkappa}) is clearly invariant under permutations of the $x_i$
thereby satisfying the Weyl group invariance hypothesis of Proposition~\ref{Pro-WG-invariance}. It follows
that on the right hand side of (\ref{eqn-gln-gf3-xkappa}) $\x^\kappa$ may be replaced by $\ch^{GL(n)}_{\kappa}(\x)$.
Furthermore, thanks to (\ref{eqn-ch-wkappa}) $\ch^{GL(n)}_{\kappa}(\x)$ may itself be replaced by $\sgn(w)\ch^{GL(n)}_{w\cdot\kappa}$,
for any permutation $w$, where $w\cdot\kappa$ denotes $w(\kappa+\rho)-\rho$ with $\rho=(n-1,\ldots,1,0)$. This dot action of the Weyl group 
is generated by that of $w_\alpha$ with $\alpha$ a simple root as given in Table~\ref{Tab-wkappa}.
In the case of $GL(n)$ it can be seen that, with the restriction of the components of $\kappa$ to $\{3,2,1,0\}$,
all possible pairs of consecutive components of $\kappa\in\Lambda$ in strictly increasing order, as opposed to the weakly
decreasing order required of $\lambda\in\Lambda^+$, transform under the dot action of $w_{\alpha_i}$ as shown below:
\begin{equation}\label{eqn-gln-wi}
\begin{array}{|l|llllll|}
\hline
(\kappa_i,\kappa_{i+1})&(0,3)&(1,3)&(2,3)&(0,2)&(1,2)&(0,1)\cr
(\kappa_{i+1}-1,\kappa_{i}+1)&(2,1)&(2,2)&(2,3)&(1,1)&(1,2)&(0,1)\cr
\hline
\end{array}
\end{equation}

Iterating these transformations and keeping track of their signature factors
one finds that in general $\lambda\in\Lambda_{GL(n)}^+$ is of the general form $(3^p,2^q,1^r,0^{n-p-q-r})$
and may only be built from subsequences $\tau= w\cdot\sigma$ with $\sigma$
certain specific elementary subsequences of $\kappa$, where it has been convenient to 
write $\tau=w\cdot\sigma$ if $w((\ldots,\sigma,\ldots)+\rho)-\rho=(\ldots,\tau,\ldots)$. 
These are displayed in Table~\ref{Tab-seq-gln} along with the signature $\sgn(w)$ and the 
corresponding contribution of $a(\sigma)$ to $a(\kappa)$. 
All other subsequences $\sigma$ of $\kappa$ are prohibited in that at some stage in the iteration process
one has to conclude that $\ch_\kappa^{GL(n)}=0$ as a result of the occurrence of one or other
of the pairs of consecutive components $(2,3)$, $(1,2)$ or $(0,1)$ since these pairs are invariant under the 
dot action of $w_\alpha$ with $\sgn(w_\alpha)=-1$.

\begin{table}[ht]
\begin{center}
\begin{tabular}{|l|l|l|}
\hline
$\sigma$&$\tau=w\cdot\sigma$&$\sgn(w)\,a(\sigma)$ \cr
\hline\hline$(3)$&$(3)$&$+a_3$\cr
\hline
$(2)$&$(2)$&$+a_2$\cr
$(1,3)$&$(2,2)$&$-a_1a_3$\cr
$(0,3,3)$&$(2,2,2)$&$+a_0a_3^2$\cr
\hline
$(2,1)$&$(2,1)$&$+a_2a_1$\cr
$(0,3)$&$(2,1)$&$-a_0a_3$\cr
\hline
$(1)$&$(1)$&$+a_1$\cr
$(0,2)$&$(1,1)$&$-a_0a_2$ \cr
$(0,0,3)$&$(1,1,1)$&$+a_0^2a_3$\cr
\hline
$(0)$&$(0)$&$+a_0$\cr
\hline
\hline
\end{tabular}
\end{center}
\medskip
\caption{Elementary subsequences $\tau=w\cdot\sigma$ of $\lambda$ along with $\sgn(w)$ and the contribution of $a(\sigma)$ to $a(\kappa)$ in the case $GL(n)$.} 
\label{Tab-seq-gln}
\end{table}

Initial subsequences of $3$'s and trailing subsequences of $0$'s are left invariant by the
re-ordering transformations, thereby giving rise to the factors $a_3^p$ and $a_0^{n-p-q-r}$ appearing on the right hand 
side of (\ref{eqn-gln-gf3}). The subsequences $\tau=(3,2)$ and $(1,0)$ of $\lambda$ may only be formed in one way, unlike $\tau=(2,1)$ 
which can arise from both $\sigma=(2,1)$ and $(0,3)$. This observation is responsible for the appearance of the
two terms in Table~\ref{Tab-seq-gln}. 
Taking into account the weighting $\sgn(w)\,a(\sigma)$ of each of the subsequences of Table~\ref{Tab-seq-gln}
and the required number of entries $3$, $2$, $1$ and $0$ in $\lambda$ one arrives at the recurrence relations
(\ref{eqn-gln-Q}) and (\ref{eqn-gln-R}), thereby completing the derivation of the expansion (\ref{eqn-gln-gf3}).
\qed

In the special cases obtained by setting $a_2=a_3=0$ or just $a_3=0$ in Theorem~\ref{The-gln} we have 
\begin{Corollary}\label{Cor-gln}
For $\a=(a_0,a_1)$ and $\a=(a_0,a_1,a_2)$
\begin{align}
 \prod_{i=1}^n (a_0+a_1x_i) &=\ds \sum_{r=0}^{n}\, a_1^r\,a_0^{n-r}\ \ch^{GL(n)}_{(1^r)}(\x)\,;\label{eqn-gln-gf1}\\.
 \prod_{i=1}^n (a_0+a_1x_i+a_2x_i^2) &=\ds  \sum_{q=0}^{n} \, \sum_{r=0}^{n-q}\ a_2^q\,\phi_{r}(\a)\,a_0^{n-q-r}\ \ch^{GL(n)}_{(2^q,1^r)}(\x)\,,\label{eqn-gln-gf2-phi}
\end{align}
where
\begin{equation}\label{eqn-gln-phi-ABC}
\phi_0(\a)=1,~~\phi_1(\a)=a_1~~\mbox{and}~~\phi_r(\a)=a_1\phi_{r-1}(\a)-a_0a_2\phi_{r-2}(\a)~~\mbox{for $r\geq2$}\,.
\end{equation}
or equivalently
\begin{equation}
\phi_r(\a)= [t^r]\, 1/(1-a_1\,t +a_0a_2\,t^2)~~\mbox{for all $r\geq0$}\,.
\end{equation}
where $[t^r]\,P(t)$ signifies the coefficient of $t^r$ in $P(t)$ expanded as a power series in $t$.
\end{Corollary}

Results obtained from the recurrence relations of Corollary~\ref{Cor-gln} and Theorem~\ref{The-gln} 
for some specific values of $\a=(1,a_1,a_2)$ and $\a=(1,a_1,a_2,a_3)$ are offered in 
Tables~\ref{Tab-gln-gf2-values} and~\ref{Tab-gln-gf3-values}, respectively. They consist
mainly of results obtained previously in~\cite{KWY,YW} and by different means in~\cite{LO2}.
In these tabulations use has been made of the following notation:
\begin{equation}
\begin{array}{rl}
(i)&F_0=1,~~F_1=1~~\mbox{and}~~F_k=F_{k-1}+F_{k-2}~~\mbox{for $k\geq2$};\cr
(ii)&G_0=0,~~G_1=0,~~G_2=1~~\mbox{and}~~G_k=G_{k-1}+G_{k-2}+G_{k-3}~~\mbox{for $k\geq3$};\cr
(iii)&H_0=0,~~H_1=0,~~H_2=1~~\mbox{and}~~H_k=-H_{k-1}-H_{k-2}+H_{k-3}~~\mbox{for $k\geq3$},\cr
\end{array}
\end{equation}
see~\cite{OEIS}~(i)~A000045, (ii)~A000073, (iii)~A57597,
where $F_k$ are the Fibonacci numbers, $G_k$ the Tribonacci numbers and $H_k$ a variation on the latter.

\begin{table}[ht]
\begin{center}
\begin{tabular}{|c||l|l|}
\hline
$\a$&$a_2^q\,\phi_{q,r}(\a)$&\cr 
\hline
\hline
$(1,0,1)$
&$\begin{array}{rl}
     (-1)^{r/2}&\mbox{if $r=0\mod 2$}\cr
			0&\mbox{if $r=1\mod 2$}\cr
   \end{array}$
&\begin{footnotesize}$\begin{array}{l}\mbox{\cite{YW} Table 2 $V^+$}\cr\mbox{\cite{LO2} (3.2)}\cr\end{array}$\end{footnotesize}\cr 
\hline
$(1,0,\ov1)$
&$\begin{array}{rl}
     (-1)^{q}&\mbox{if $r=0\mod 2$}\cr
			0&\mbox{if $r=1\mod 2$}\cr
   \end{array}$
&\begin{footnotesize}$\begin{array}{l}\mbox{\cite{YW} Table 1 $V$}\cr\mbox{\cite{LO2} (3.6)}\cr\end{array}$\end{footnotesize}\cr
\hline
$(1,1,1)$
&$\begin{array}{rl}
                      1&\mbox{if $r=0,1\mod 6$}\cr
									 		0&\mbox{if $r=2,5\mod 6$}\cr
										 -1&\mbox{if $r=3,4\mod 6$}\cr
                        \end{array}$
&~\begin{footnotesize}\cite{YW}~(33)\end{footnotesize}\cr
\hline
$(1,\ov1,1)$
&$\begin{array}{rl}
                      1&\mbox{if $r=0\mod 3$}\cr
									  	0&\mbox{if $r=2\mod 3$}\cr
										 -1&\mbox{if $r=1\mod 3$}\cr
                        \end{array}$
&~\begin{footnotesize}\cite{YW}~(34)\end{footnotesize}\cr
\hline
$(1,1,\ov1)$
&$(-1)^q\,F_{r+1}$
&~\begin{footnotesize}\cite{KWY}~(5.23)\end{footnotesize}\cr
\hline
$(1,\ov1,\ov1)$&$(-1)^{q+r}\,F_{r+1}$&\cr
\hline
$(1,2,1)$
&$r+1$
&~\begin{footnotesize}\cite{LO2}~(3.11)\end{footnotesize}\cr
\hline
$(1,\sqrt{2},1)$&$\begin{array}{rl}
                     1&\mbox{if $r=0,2\mod8$}\cr
										 \sqrt{2}&\mbox{if $r=1\mod8$}\cr
										  0&\mbox{if $r=3,7\mod8$}\cr
											-1&\mbox{if $r=4,6\mod8$}\cr
										  -\sqrt{2}&\mbox{if $r=5\mod8$}\cr
											\end{array}$
									 &\cr
\hline
$(1,3,1)$&$F_{2r+2}$&\cr
\hline
\end{tabular}
\end{center}
\medskip
\caption{The $GL(n)$ coefficients $\phi_{r}(\a)$ in (\ref{eqn-gln-gf2-phi}) for various $\a=(1,a_1,a_2)$} 
\label{Tab-gln-gf2-values}
\end{table}

\begin{table}[ht]
\begin{center}
\begin{tabular}{|c||l|l|}
\hline
$\a$&$\psi_{q,r}(\a)$&\cr 
\hline
\hline
$(1,0,0,1)$&$\begin{array}{rl} 1&\mbox{if $q=0\mod 3$ and $r=0\mod3$}\cr
                           -1&\mbox{if $q=1\mod 3$ and $r=1\mod3$}\cr
													   0&\mbox{otherwise}\cr
														\end{array}$
&\mbox{\begin{footnotesize}$\begin{array}{c}\mbox{\cite{YW}~(29a)}\cr\end{array}$\end{footnotesize}}\cr
\hline
$(1,0,0,\ov1)$&$\begin{array}{rl} (-1)^{(3p+r)/3}&\mbox{if $q=0\mod 3$ and $r=0\mod3$}\cr
                            (-1)^{(3p+r-1)/3}&\mbox{if $q=1\mod 3$ and $r=1\mod3$}\cr
														0&\mbox{otherwise}\cr
														\end{array}$
&\begin{footnotesize}$\begin{array}{c}\mbox{\cite{YW}~(29b)}\cr \end{array}$\end{footnotesize}\cr
\hline
$(1,1,1,1)$&
$\begin{array}{l}
       \psi_{q+4,r}=\psi_{q,r+4}=\psi_{q,r}~~\mbox{with}\cr
			 \psi_{q,r}=\mbox{\begin{footnotesize}$\begin{array}{|c||r|r|r|r|}
                      \hline
                      q\backslash r&0&1&2&3\cr
                      \hline\hline                      
											0&1&1&0&0\cr
                       \hline
											1&1&0&-1&0\cr
                       \hline
                      2&0&-1&-1&0\cr
                       \hline
                      3&0&0&0&0\cr
                       \hline
                      \end{array}$\end{footnotesize}}\cr\cr
\end{array}$
&\begin{footnotesize}$\begin{array}{l}\mbox{\cite{KWY}~(5.15)}\cr\mbox{\cite{LO2}~(4.8)}\cr\end{array}$\end{footnotesize}\cr
\hline
$(1,\ov1,\ov1,1)$&
$\begin{array}{l}
       \psi_{q+2,r}=\psi_{q,r+2}=\psi_{q,r}~~\mbox{with}\cr
			 \psi_{q,r}=\mbox{\begin{footnotesize}$\begin{array}{|c||l|l|}
                      \hline
                      q\backslash r&0&1\cr
                      \hline\hline                      
											0&(q+r+2)/2&-(r+1)/2\cr
                      \hline
											1&-(q+1)/2&0\cr
                      \hline
                      \end{array}$\end{footnotesize}}\cr\cr
\end{array}$
&\begin{footnotesize}$\begin{array}{l}\mbox{\cite{KWY}~(5.21)}\cr
                                     \mbox{see also}\cr
															        \mbox{\cite{LO2}~(4.18)}\cr\end{array}$\end{footnotesize}\cr
\hline
$(1,\ov1,1,1)$&
$G_{q+2}H_{r+2}-G_{q+1}H_{r+1}$&\cr
\hline
$(1,1,\ov1,1)$&
$H_{q+2}G_{r+2}-H_{q+1}G_{r+1}$&\cr
\hline
$(1,2,2,1)$&
$\begin{array}{l}
       \psi_{q+6,r}=\psi_{q,r+6}=\psi_{q,r}~~\mbox{with}\cr\cr
			 \psi_{q,r}=\mbox{\begin{footnotesize}$\begin{array}{|c||r|r|r|r|r|r|}
                      \hline
                      q\backslash r&0&1&2&3&4&5\cr
                      \hline
											\hline                      
											0&1&2&2&1&0&0\cr
                      1&2&3&2&0&-1&0\cr
                      2&2&2&0&-2&-2&0\cr
                      3&1&0&-2&-3&-2&0\cr
                      4&0&-1&-2&-2&-1&0\cr
                      5&0&0&2&1&0&0\cr
                      \hline
                      \end{array}$\end{footnotesize}}\cr\cr
\end{array}$&\cr
\hline
$(1,3,3,1)$&
$(q+r+2)(q+1)(r+1)/2$&\cr
\hline
\end{tabular}
\end{center}
\medskip
\caption{The $GL(n)$ coefficients $\psi_{q,r}(\a)$ in (\ref{eqn-gln-gf3}) for various $\a=(1,a_1,a_2,a_3)$} 
\label{Tab-gln-gf3-values}
\end{table}

In principle all of these results may be arrived at by exploiting the dual Cauchy identity~\cite{Mac,BG} given 
by (\ref{eqn-gln-glm}). In the case $m=2$ and $a_2=1$ one has
\begin{equation}
 \prod_{i=1}^n (a_0+a_1x_i+x_i^2)=\prod_{i=1}^n \prod_{j=1}^2 (x_i+y_j)=\sum_{\lambda\in(2^n)}\ \ch_\lambda^{GL(n)}(\x)\ \ch_{\tilde\lambda}^{GL(2)}(y_1,y_2) 
\end{equation}
with $a_0=y_1y_2$ and $a_1=(y_1+y_2)$. Clearly $\lambda$ is necessarily of the form $(2^q,1^r)$ with $\tilde\lambda=(s+r,s)$ with $s=n-q-r$.
In this case $\ch_{(s+r,s)}^{GL(2)}(y_1,y_2)=(y_1y_2)^s\ch_{(r)}^{GL(2)}(y_1,y_2)=a_0^{n-q-r}\phi_r(\a)$, so that in the notation of (\ref{eqn-gln-gf2-phi})
with $a_2=1$
\begin{equation}
     \phi_r(\a)=\ch_{(r)}^{GL(2)}(y_1,y_2)=\sum_{j=0}^r y_1^{r-j}\,y_2^j\,.
\end{equation}
Setting $y_1=e^{i\pi/k}$ and $y_2=e^{-i\pi/k}$ with $k=1,2,3,4$ and $6$ yields the results of Table~\ref{Tab-gln-gf2-values} 
in the cases $\a=(1,a_1,1)$ with $a_1=-2,0,1,\sqrt{2}$ and $\sqrt{3}$, while the even simpler case $y_1=y_2=1$ corresponds to $a_1=2$.

Similarly, in the case $m=3$ and $a_3=1$ one has
\begin{equation}
 \prod_{i=1}^n (a_0+a_1x_i+a_2x_i^2+x_i^3)=\prod_{i=1}^n \prod_{j=1}^3 (x_i+y_j)=\sum_{\lambda\in(3^n)}\ \ch_\lambda^{GL(n)}(\x)\ \ch_{\tilde\lambda}^{GL(3)}(y_1,y_2,y_3) 
\end{equation}
with $a_0=y_1y_2y_3$, $a_1=y_1y_2+y_1y_3+y_2y_3$ and $a_1=y_1+y_2+y_3$. This time $\lambda$ is of the form $(3^p,2^q,1^r)$ with $\tilde\lambda=(s+q+r,s+r,s)$ with $s=n-p-q-r$.
In this case $\ch_{(s+q+r,s+r,s)}^{GL(3)}(y_1,y_2,y_3)=(y_1y_2y_3)^s\ch_{(q+r,r)}^{GL(3)}(y_1,y_2,y_3)=a_0^{n-q-r}\psi_{q,r}(\a)$, so that in the notation of (\ref{eqn-gln-gf3}) with $a_3=1$
\begin{equation}
     \psi_{q,r}(\a)=\ch_{(q+r,r)}^{GL(3)}(y_1,y_2,y_3)\,.
\end{equation}
In general it is not so easy to evaluate this for all $(q,r)$. In particular it does not lend itself well to the evaluation of
$\psi_{q,r}(\a)$ in the cases $\a=(1,-1,1,1)$ and $(1,1,-1,1)$ of Table~\ref{Tab-gln-gf3-values}. 
However, in the periodic cases obtained by setting $(y_1,y_2,y_3)=(1,e^{i2\pi/k},e^{-i2\pi/k})$
with $k=2,3,4$ and $6$ one recovers the results of Table~\ref{Tab-gln-gf3-values} in the cases $\a=(1,a,a,1)$ with $a=-1,0,1$ and $2$, respectively.
as well as the simpler case $y_1=y_2=y_3$ corresponding to $a=3$, in which case $\psi_{q,r}(\a)$ is just the dimension of the $GL(3)$ representation of highest weight $(q+r,r)$.

\section{Symplectic group character identities}\label{sec-spn}
Let $n\in\N$ be fixed, and let $\x=(x_1,x_2,\ldots,x_n)$ and $\ov\x=(\ov{x}_1,\ov{x}_2,\ldots,\ov{x}_n)$ with 
$\ov{x}_i=x_i^{-1}$ for $i=1,2,\ldots,n$, and let $\ch^{Sp(2n)}_\lambda(\x,\ov\x)$ denote the character 
of the irreducible representation of $Sp(2n)$ of highest weight $\lambda$ evaluated on a group element with eigenvalues $(\x,\ov\x)$.
Then we have the following 

\begin{Theorem}\label{The-spn}
For all $\a=(a_0,a_1,a_2)$ we have
\begin{equation}\label{eqn-spn-gf}
   \prod_{i=1}^n (a_0+a_1(x_i+\ov{x}_i)+a_2(x_i^2+\ov{x}_i^2))
	     =\sum_{p=0}^n \sum_{q=0}^{n-p} \delta_{r,n-p-q}\ a_2^p\ \psi_{q,r}(\a)\ \ch^{Sp(2n)}_{(2^p,1^q,0^r)}(\x,\ov\x)\,,
\end{equation}
where
\begin{equation}\label{eqn-spn-psi}
  \psi_{q,r}(\a)=\chi_{q,r}(\a)-\chi_{q,r-1}(\a)\,a_2\,,
\end{equation}
with
\begin{equation}\label{eqn-spn-chi}
\chi_{q,r}(\a)= Q_q\,R_r+a_2^3\,Q_{q-2}\,R_{r-1}+a_1\,Q_{q-1}\sum_{s=1}^r (-a_2)^s\,R_{r-s}\,,
\end{equation}
where~~$Q_q=0$ if $q<0$,~~$Q_0=1$ and 
\begin{equation}\label{eqn-spn-Q}
Q_q=a_1\,Q_{q-1}-a_0a_2\,Q_{q-2}+a_1a_2^2\,Q_{q-3}-a_2^4\,Q_{q-4}~~\mbox{if $q\geq1$}\,,
\end{equation}
while~~$R_r=0$ if $r<0$,~~$R_0=1$ and
\begin{equation}\label{eqn-spn-R}
R_r=a_0\,R_{r-1}-a_1^2\,R_{r-2}-a_0a_2^2\,R_{r-3}+a_2^4\,R_{r-4}   
\ds-2a_1^2\sum_{s=1}^{r-2}(-a_2)^s\,R_{r-s-2}~~\mbox{if $r\geq1$}\,.   
\end{equation}
\end{Theorem}

\noindent{\bf Proof}:~~
\begin{equation}\label{eqn-spn-gf-xkappa}
   \prod_{i=1}^n (a_0+a_1(x_i+\ov{x}_i)+a_2(x_i^2+\ov{x}_i^2))=\sum_{\kappa}  a_2^k\,a_1^\ell\,a_0^{n-k-\ell} \x^\kappa\,,
\end{equation}
where $\kappa_j\in\{\ov2,\ov1,0,1,2\}$ for $j=1,2,\ldots,n$ and
\begin{equation}
\begin{array}{rcl}
     k&=&\#\{\kappa_j\in\{2,\ov2\}|j=1,2,\ldots,n\};\cr
		 \ell&=&\#\{\kappa_j\in\{1,\ov1\}|j=1,2,\ldots,n\}.\cr
\end{array}
\end{equation}
The left hand side of (\ref{eqn-spn-gf-xkappa}) is clearly invariant under permutations of the $x_i$ and sign changes of their components,
thereby satisfying the Weyl group invariance hypothesis of Proposition~\ref{Pro-WG-invariance}. It follows from (\ref{eqn-xkappa-chkappa})
that on the right hand side of (\ref{eqn-spn-gf-xkappa}) $\x^\kappa$ may be replaced by $\ch^{Sp(2n)}_{\kappa}(\x,\ov\x)$.

Referring to the case $Sp(2n)$ of Table~\ref{Tab-wkappa} it can then be seen  
that pairs of consecutive components of $\kappa$ in non-standard order transform under the dot action of 
$w_{\alpha_i}$ as shown below for $i=1,2,\ldots,n-1$:
\begin{equation}\label{eqn-spn-wi}
\begin{array}{|l|llllll|}
\hline
(\kappa_i,\kappa_{i+1})&(\ov2,0)&(\ov2,1)&(\ov2,2)&(\ov1,1)&(\ov1,2)&(0,2)\cr
(\kappa_{i+1}-1,\kappa_{i}+1)&(\ov1,\ov1)&(0,\ov1)&(1,\ov1)&(0,0)&(1,0)&(1,1)\cr
\hline
(\kappa_i,\kappa_{i+1})&(\ov2,\ov1)&(\ov1,0)&(0,1)&(1,2)&&\cr
(\kappa_{i+1}-1,\kappa_{i}+1)&(\ov2,\ov1)&(\ov1,0)&(0,1)&(1,2)&&\cr
\hline
\end{array}
\end{equation}
while the dot action of $w_{\alpha_n}$ with $\alpha=2\epsilon_n$ transforms the $n$th component of $\kappa$ as indicated below: 
\begin{equation}\label{eqn-spn-wn}
\begin{array}{|l|ll|}
\hline
(\ldots,\kappa_n)&(\ldots,\ov2)&(\ldots,\ov1)\cr
(\ldots,-\kappa_n-2)&(\ldots,0)&(\ldots,\ov1)\cr
\hline
\end{array}
\end{equation} 

The transformations (\ref{eqn-spn-wi}) and (\ref{eqn-spn-wn}) lead inexorably to the transformations and weightings of elementary subsequences 
given in Table~\ref{Tab-seq-spn} in which $t$ may be any positive integer.
\begin{table}[ht]
\begin{center}
\begin{tabular}{|l|l|l|}
\hline
$\sigma$&$\tau=w\cdot\sigma$&$\sgn(w)\,a(\sigma)$ \cr
\hline\hline$(2)$&$(2)$&$+a_2$\cr
\hline
$(2,1)$&$(2,1)$&$+a_1a_2$\cr
\hline
$(1)$&$(1)$&$+a_1$\cr
$(0,2)$&$(1,1)$&$-a_0a_2$ \cr
$(\ov1,2,2)$&$(1,1,1)$&$+a_1a_2^2$\cr
$(\ov2,2,2,2)$&$(1,1,1,1)$&$-a_2^4$\cr
\hline
$(1,0)$&$(1,0)$&$+a_0a_1$\cr
$(\ov1,2)$&$(1,0)$&$-a_1a_2$\cr
$(\ov2,2,2)$&$(1,1,0)$&$+a_2^3$\cr
\hline
$((\ov2,2)^t,1)$&$(1,0^{2t})$&$+a_1a_2^{2t}$\cr
$((\ov2,2)^t,\ov1,2)$&$(1,0^{2t+1})$&$-a_1a_2^{2t+1}$\cr
\hline
$(0)$&$(0)$&$+a_0$\cr
$(\ov1,1)$&$(0,0))$&$-a_1^2$\cr
$(\ov2,0,2)$&$(0,0,0))$&$-a_0a_2^2$\cr
$(\ov2,\ov2,2,2)$&$(0,0,0,0))$&$+a_2^4$\cr
\hline
$(\ov1,(\ov2,2)^t,\ov1,2)$&$(0^{2t+3})$&$+a_1^2a_2^{2t+1}$\cr
$(\ov2,1,(\ov2,2)^t,1)$&$(0^{2t+3})$&$+a_1^2a_2^{2t+1}$\cr
$(\ov1,(\ov2,2)^{t+1},1)$&$(0^{2t+4})$&$-a_1^2a_2^{2t+2}$\cr
$(\ov2,1,(\ov2,2)^t,\ov1,2)$&$(0^{2t+4})$&$-a_1^2a_2^{2t+2}$\cr
\hline
\hline
\end{tabular}
\end{center}
\medskip
\caption{Elementary subsequences $\tau=w\cdot\sigma$ of $\lambda$ along with $\sgn(w)$ and the contribution of 
$a(\sigma)$ to $a(\kappa)$ in the case $Sp(2n)$.} 
\label{Tab-seq-spn}
\end{table}

In deriving the $t$-dependent transformations of Table~\ref{Tab-seq-spn} from (\ref{eqn-spn-wi})
it should be noted that $(\ov2,2)^t$ maps to $(1,\ov1)^t$ with signature factor $(-1)^t$,
while $(\ov1,1)^t$ maps to $(0,0)^t$ with the same signature factor $(-1)^t$.

Any initial sequence of $2$'s in $\kappa$ is left invariant through the dot action of $w_{\alpha_i}$ for all $i$. 
This is the origin of the factor $a_2^p$ on the right hand side of (\ref{eqn-spn-gf}). 

The transformations of (\ref{eqn-spn-wn}) imply that 
that in (\ref{eqn-spn-gf-xkappa}) all sequences $\kappa$ ending in $\ov1$ can be ignored, 
while any sequence $\kappa$ ending in $\ov2$ gives a contribution obtained by changing the sign of that
arising if $\ov2$ is replaced by $0$. This observation has been taken into account through the
inclusion of the second term in the definition~(\ref{eqn-spn-psi}) of $\psi_{q,r}(\a)$.

The terms with $\tau$ of the form $(1^\ell,0^m)$ in Table~\ref{Tab-seq-spn} give rise to the expression
(\ref{eqn-spn-chi}) for $\chi_{q,r}(\a)$ in terms of $Q_{q-\ell}$ and $R_{r-m}$, while the expressions for the latter
are derived from the entries with $\tau$ of the form $(1^\ell)$ and $(0^m)$, respectively in Table~\ref{Tab-seq-spn}.
These observations complete the derivation of (\ref{eqn-spn-gf}).
\qed 

As special cases of Theorem~\ref{The-spn} obtained by setting $a_2=0$ and $a_1=0$ we have 
\begin{Corollary}\label{Cor-spn-AB-AC}
For $\a=(a_0,a_1)$ and $\a=(a_0,0,a_2)$ we have
\begin{align}
 \ds \prod_{i=1}^n (a_0+a_1(x_i+\ov{x}_i))&=\ds \sum_{q=0}^{n}\,\delta_{r,n-q}\,a_1^q\,\phi_r(\a)\,\ch^{Sp(2n)}_{(1^q,0^r)}(\x,\ov\x)\,;\label{eqn-spn-AB}\\
 \ds \prod_{i=1}^n (a_0+a_2(x_i^2+\ov{x}_i^2))&=\ds \sum_{p=0}^n \sum_{q=0}^{n-p}\,\delta_{r,n-p-q}\,a_2^p\,\psi_{q,r}(\a)\,\ch^{Sp(2n)}_{(2^p,1^q,0^r)}(\x,\ov\x)\,,\label{eqn-spn-AC}
\end{align}
with
\begin{equation}\label{eqn-spn-phi-AB}
\phi_0(\a)=1,\phi_1(\a)=a_0~~\mbox{and}~~\phi_r(\a)=a_0\,\phi_{r-1}(\a)-a_1^2\,\phi_{r-2}(\a)\,,
\end{equation}
and
\begin{equation}\label{eqn-spn-psi-AC}
\begin{array}{l}
\psi_{q,r}(\a)=0~~\mbox{if $q<0$ or $q=1\mod2$ or $r<0$}\,; \cr
\psi_{0,0}(\a)=1~\mbox{and}~~\psi_{0,1}(\a)=a_0-a_2\,;\cr  
\psi_{0,r}(\a)=a_0\psi_{0,r-1}(\a)-a_0a_2^2\psi_{0,r-3}(\a)+a_2^4\psi_{0,r-4}(\a)~~\mbox{for $r\geq 2$}\,; \cr
\psi_{2,r}(\a)=-a_0a_2\psi_{0,r}(\a)+a_2^3\psi_{0,r-1}~\mbox{for $r\geq 0$}\,; \cr
\psi_{q,r}(\a)=-a_0a_2\psi_{q-2,r}(\a)-a_2^4\psi_{q-4,r}(\a)~~\mbox{for $q=0\mod2$ with $q\geq 4$ and $r\geq 0$}\,. \cr
\end{array}
\end{equation}
\end{Corollary}

\noindent{\bf Proof}:
The first case corresponds to setting $a_2=0$ in (\ref{eqn-spn-gf}). This immediately gives $\psi_{q,r}(\a)=Q_qR_r$
with $Q_q=a_1^q$ and $R_r=a_0R_{r-1}-a_1^2R_{r-2}$, so that $\psi_{q,r}=a_1^q\phi_{q}(\a)$ with 
$\phi_0(\a)=1$ and $\phi_r(\a)=a_0\,\phi_{r-1}(\a)-a_1^2\,\phi_{r-2}(\a)$, as required.

In the second case, one sets $a_1=0$ in (\ref{eqn-spn-gf}). This gives (i) $\psi_{q,r}(\a)=\chi_{q,r}(\a)-a_2\chi_{q,r-1}(\a)$; 
(ii) $\chi_{q,r}(\a)=Q_qR_r+a_2^3Q_{q-2}R_{r-1}$; (iii) $Q_q=-a_0a_2Q_{q-2}-a_2^4Q_{q-4}$ for $q\geq1$; (iv) $R_r=a_0R_{r-1}-a_0a_2^2R_{r-3}+a_2^4R_{r-4}$ 
for $r\geq1$.
Since $Q_q=0$ for $q<0$ and $R_r=0$ for $r<0$ it follows from (i) and (ii) that $\psi_{q,r}(\a)=0$ if either $q<0$ or $r<0$.
From (iii) and the fact that $Q_1=a_1=0$ it is clear $Q_q=0$ if $q$ is odd. Then (i) and (ii) imply that $\psi_{q,r}(\a)=0$ if $q$ is odd.
The conditions $Q_0=R_0=1$ are then sufficient to ensure that $\phi_{0,r}(\a)=R_r$ and $\psi_{0,r}(\a)=R_r-a_2R_{r-1}$. 
The recurrence relation (iv) for $R_r$ then leads directly to the required expression for $\psi_{0,r}(\a)$ for all $r\geq0$.
The fact that $Q_2=-a_0a_2$ and $Q_0=1$ implies that $\chi_{2,r}(\a)=-a_0a_2R_r+a_2^3R_{r-1}$ so that from (i) we have
$\psi_{2,r}(\a)=-a_0a_2R_r+a_2^3R_{r-1}+a_0a_2^2R_{r-1}-a_2^4R_{r-2}=-a_0a_2\psi_{0,r}(\a)+a_2^3\psi_{0,r-1}(\a)$, again as required.
Finally, from (i) and (ii) we find $\psi_{q,r}(\a)+a_0a_2\psi_{q-2,r}(\a)+a_2^4\psi_{q-4,r}(\a)=Z_qR_r+a_2^3Z_{q-2}R_{r-1}-a_2Z_qR_{r-1}-a_2^4Z_{q-2}R_{r-2}$ 
where $Z_q=Q_q+a_0a_2Q_{q-2}+a_2^4Q_{q-4}=0$ for all $q\geq1$ by virtue of (iii). This yields the required recurrence relation for $\psi_{q,r}(\a)$ 
for all $r\geq0$ and $q\geq3$, that is for $q\geq4$ since we only require the case $q=0\mod2$.
\qed

Some results obtained from Theorem~\ref{The-spn} and Corollary~\ref{Cor-spn-AB-AC} for various specific values 
of $\a=(a_0,1)$, $(a_0,0,1)$ and $(a_0,1,1)$ are offered in Tables~\ref{Tab-phi-values}, \ref{Tab-spn-psi-A01} and \ref{Tab-spn-psi-A11},
respectively.


\section{Odd orthogonal group character identities}\label{sec-oon}
For any partition $\lambda$ let $\ch^{SO(2n+1)}_\lambda(\x,\ov\x,1)$ denote the character 
of the irreducible representation of $SO(2n+1)$ of highest weight $\lambda$ evaluated on a group element with eigenvalues $(\x,\ov\x,1)$.
Then we have

\begin{Theorem}\label{The-oon}
For all $\a=(a_0,a_1,a_2)$ 
\begin{equation}\label{eqn-oon-gf}
   \prod_{i=1}^n (a_0+a_1(x_i+\ov{x}_i)+a_2(x_i^2+\ov{x}_i^2))
	     =\sum_{p=0}^n \sum_{q=0}^{n-p} \delta_{r,n-p-q}\ a_2^p\ \psi_{q,r}(\a)\ \ch^{SO(2n+1)}_{(2^p,1^q,0^r)}(\x,\ov\x,1)\,,
\end{equation}
where
\begin{equation}\label{eqn-oon-psi}
      \psi_{q,r}(\a)= \sum_{s=0}^r \chi_{q,r-s}(\a)\,S_s\ +\ \chi_{q-1,0}(\a)\,(-a_2)^{r+1} \,,
\end{equation}
with $\chi_{q,r}(\a)$ defined as in (\ref{eqn-spn-chi}) in terms of $Q_q$ and $R_r$ as given by (\ref{eqn-spn-Q}) and (\ref{eqn-spn-R}), respectively, and 
\begin{equation}\label{eqn-oon-S}
    S_0=1,~~S_1=-a_1,~~S_2=2a_1a_2-a_2^2,~~\mbox{and $S_s=-2a_1(-a_2)^{s-1}$ if $s\geq3$}\,.
\end{equation}
\end{Theorem}

\noindent{\bf Proof}:~~
As in the symplectic case we have
\begin{equation}\label{eqn-oon-gf-xkappa}
   \prod_{i=1}^n (a_0+a_1(x_i+\ov{x}_i)+a_2(x_i^2+\ov{x}^2))=\sum_{\kappa}  a_2^k\,a_1^\ell\,a_0^{n-k-\ell} \x^\kappa\,,
\end{equation}
where $\kappa_j\in\{\ov2,\ov1,0,1,2\}$ for $j=1,2,\ldots,n$ and
\begin{equation}
\begin{array}{rcl}
     k&=&\#\{\kappa_j\in\{2,\ov2\}|j=1,2,\ldots,n\};\cr
	\ell&=&\#\{\kappa_j\in\{1,\ov1\}|j=1,2,\ldots,n\}.\cr
\end{array}
\end{equation}
Since $SO(2n+1)$ shares the same Weyl group as $Sp(2n)$ the left hand side of (\ref{eqn-oon-gf-xkappa}) 
again satisfies the Weyl group invariance hypothesis of Proposition~\ref{Pro-WG-invariance}. It follows
that on the right hand side $\x^\kappa$ may be replaced by $\ch^{SO(2n+1)}_{\kappa}(\x,\ov\x,1)$.
Moreover the tansformations of (\ref{eqn-spn-wi}) and their iterated consequences given in Table~\ref{Tab-seq-spn}
still apply. 

However, what is different in the $SO(2n+1)$ case is the action of $w_{\alpha_n}$. 
In this case $\alpha_n=\epsilon_n$ and the dot action of $w_{\epsilon_n}$  transforms the 
$n$th component of $\kappa$ as shown while leaving all other components unchanged. 
\begin{equation}\label{eqn-oon-wn}
\begin{array}{|l|ll|}
\hline
(\ldots,\kappa_n)&(\ldots,\ov2)&(\ldots,\ov1)\cr
(\ldots,-\kappa_n-1)&(\ldots,1)&(\ldots,0)\cr
\hline
\end{array}
\end{equation}
This implies that 
all sequences $\kappa$ ending in $\ov2$ or $\ov1$ can be replaced by sequences ending in $1$ or $0$, respectively, 
while retaining an additional signature factor $-1$. 

Combining these observations about the final or terminating entry in any sequence $\kappa$ with those of (\ref{eqn-spn-wi}) one arrives 
at the list of transformations of terminal subsequences given in Table~\ref{Tab-seq-oon} in which $t$ may be any positive integer. 
\begin{table}[ht]
\begin{center}
\begin{tabular}{|l|l|l|}
\hline
$\sigma$&$\tau=w\cdot\sigma$&$\sgn(w)\,a(\sigma)$ \cr
\hline
\hline
$(\cdots,\ov2)$&$(\cdots,1)$&$-a_2$\cr
$(\cdots,\ov2,2)$&$(\cdots,1,0)$&$+a_2^2$\cr
\hline
$(\cdots,(\ov2,2)^{t},\ov2)$&$(\cdots,1,0^{2t})$&$-a_2^{2t+1}$\cr
$(\cdots,(\ov2,2)^{t+1})$&$(\cdots,1,0^{2t+1})$&$+a_2^{2t+2}$\cr
\hline
\hline
$(\cdots,\ov1)$&$(\cdots,0)$&$-a_1$\cr
\hline
$(\cdots,\ov1,\ov2)$&$(\cdots,0,0)$&$+a_1a_2$\cr
$(\cdots,\ov2,1)$&$(\cdots,0,0)$&$+a_1a_2$\cr
$(\cdots,\ov2,\ov2)$&$(\cdots,0,0)$&$-a_2^2$\cr
$(\cdots,\ov2,1,\ov2)$&$(\cdots,0,0,0)$&$-a_1a_2^2$\cr
\hline
$(\cdots,\ov1,(\ov2,2)^t)$&$(\cdots,0^{2t+1})$&$-a_1a_2^{2t}$\cr
$(\cdots,\ov1,(\ov2,2)^t,\ov2)$&$(\cdots,0^{2t+2)})$&$+a_1a_2^{2t+1}$\cr
$(\cdots,\ov2,1,(\ov2,2)^t)$&$(\cdots,0^{2t+2)})$&$+a_1a_2^{2t+1}$\cr
$(\cdots,\ov2,1,(\ov2,2)^t,\ov2)$&$(\cdots,0^{2t+3})$&$-a_1a_2^{2t+2}$\cr
\hline
\hline
\end{tabular}
\end{center}
\medskip
\caption{Terminating subsequences $\tau=w\cdot\sigma$ of $\lambda$ along with $\sgn(w)$ and the contribution of $a(\sigma)$ to $a(\kappa)$ 
in the case $SO(2n+1)$.} 
\label{Tab-seq-oon}
\end{table}
Here in the upper part of the table $(\cdots,)$ indicates any initial subsequence that by dint of Weyl transformations can be written in the 
form $(2^p,1^q,0^s)$, while in the lower part it must be of the form $(2^p,1^{q-1})$ with no trailing $0$'s. 
Successive transformations including signature factors are exemplified in a $t$-independent case by
$(\cdots,\ov2,1,\ov2)\mapsto-(\cdots,0,\ov1,\ov2)\mapsto(\cdots,0,\ov1,1)\mapsto-(\cdots,0,0,0)$, 
where the second step is given in Table~\ref{Tab-seq-oon} and the first and third appear in (\ref{eqn-spn-wi}). 
As in the symplectic case it is the fact that $(\ov2,2)^t$ maps to $(1,\ov1)^t$ with signature factor $(-1)^t$,
while $(\ov1,1)^t$ maps to $(0,0)^t$ with the same signature factor $(-1)^t$, that is crucial. For example,
$(\cdots,(\ov2,2)^t,\ov2)\mapsto(-1)^t(\cdots,(1,\ov1)^t,\ov2)=(-1)^{t+1}(\cdots,(1,\ov1)^t,1)
\mapsto(-1)^{t+1}(\cdots,1,(\ov1,1)^t)\mapsto-(\cdots,1,(0,0)^t)=-(\cdots,1,0^{2t})$.

Just as the transformations of Table~\ref{Tab-seq-spn} imply the validity of the recurrence relations for $Q_q$, $R_r$ and $\chi_{q,r}(\a)$,
so those of Table~\ref{Tab-seq-oon} imply the necessity of defining $\psi_{q,r}(\a)$ as in (\ref{eqn-oon-psi}) in order to encompass
all those contributions leading to $\lambda$ of the form $(\cdots,1,0^r)$ and $(\cdots,0^s)$ with the $\cdots$ signifying all terms
enumerated by $\chi_{q-1,0}(\a)=Q_{q-1}$ and $\chi_{q,r-s}(\a)$, respectively.
\qed

As special cases of Theorem~\ref{The-oon} we have

\begin{Corollary}\label{Cor-oon-AB-AC}
\begin{align}
\ds \prod_{i=1}^n (a_0+a_1(x_i+\ov{x}_i))&=\ds \sum_{q=0}^{n}\,\delta_{r,n-q}\,a_1^q\,\phi_r(\a)\,\ch^{SO(2n+1)}_{(1^q,0^r)}(\x,\ov\x,1)\,;\label{eqn-oon-AB}\\
\ds \prod_{i=1}^n (a_0+a_2(x_i^2+\ov{x}_i^2))&=\ds \sum_{p=0}^n \sum_{q=0}^{n-p}\,\delta_{r,n-p-q}\,a_2^p\,\psi_{q,r}(\a)\,\ch^{SO(2n+1)}_{(2^p,1^q,0^r)}(\x,\ov\x,1)\,,\label{eqn-oon-AC}
\end{align}
with
\begin{equation}\label{eqn-oon-phi-AB}
\begin{array}{l}
\phi_0(\a)=1,~~\phi_1=a_0-a_1~~\mbox{and}~~\phi_r(\a)=a_0\,\phi_{r-1}(\a)-a_1^2\,\phi_{r-2}(\a)~~\mbox{for $r\geq2$}\,,\cr
\end{array}
\end{equation}
and
\begin{equation}\label{eqn-oon-psi-AC}
\begin{array}{l}
\psi_{0,0}(\a)=1,~~\psi_{0,1}(\a)=a_0~~\mbox{and}~~\psi_{0,r}(\a)=a_0\psi_{0,r-1}(\a)-a_2^2\psi_{0,r-2}(\a)~~\mbox{for $r\geq2$}\,;\cr\cr  
\psi_{q,r}(\a)=\begin{cases}
                  (-a_2)^{q/2}\,\psi_{0,(q+2r)/2}(\a)&~~\mbox{for all $q,r\geq0$ if $q=0 \mod 2$}\,;\cr
									(-a_2)^{(q+1+2r)/2}\,\psi_{0,(q-1)/2}(\a)&~~\mbox{for all $q,r\geq0$ if $q=1 \mod 2$}\,.\cr
							 \end{cases}
\end{array}
\end{equation}
\end{Corollary}

\noindent{\bf Proof}: Throughout the proof it should be recalled that $Q_q=0$ if $q<0$, $R_r=0$ if $r<0$, $Q_0=R_0=S_0=1$ and $S_1=-a_1$.

{\bf Case $a_2=0$}: Setting $a_2=0$ in (\ref{eqn-oon-gf}) leaves only the terms with $p=0$ and gives $R_r=a_0R_{r-1}-a_1^2R_{r-2}$ for $r\geq1$,
$Q_q=a_1Q_{q-1}=a_1^q$ for $q\geq1$ and $\chi_{q,r}(\a)=Q_qR_r$ in (\ref{eqn-spn-R}), (\ref{eqn-spn-Q}) and (\ref{eqn-spn-chi}), respectively.
In addition, (\ref{eqn-oon-S}) implies that $S_s=0$ for $s\geq2$, so that 
(\ref{eqn-oon-psi}) gives $\psi_{q,r}(\a)=Q_qR_r-a_1Q_qR_{r-1}=a_1^q\phi_r(\a)$ with $\phi_r(\a)=R_r-a_1R_{r-1}$ for all $q,r\geq0$.
In particular $\phi_0(\a)=1$ and $\phi_1(\a)=a_0-a_1$. For $r\geq2$ the recurrence relations for $R_r$ and $R_{r-1}$ ensure 
that $\phi_r(\a)=a_0\,\phi_{r-1}(\a)-a_1^2\,\phi_{r-2}(\a)$, as required to complete the proof of (\ref{eqn-oon-phi-AB}).

{\bf Case $a_1=0$}: Setting $a_1=0$ in (\ref{eqn-oon-S}) leaves $S_0=1$, $S_1=0$ and $S_2=-a_2^2$, with $S_s=0$ for all $s\geq3$.
As a result (\ref{eqn-oon-psi}) gives 
(i) $\psi_{q,r}(\a)=(-a_2)^{r+1}\chi_{q-1,0}(\a)+\chi_{q,r}(\a)-a_2^2\chi_{q,r-2}(\a)$ for $q,r\geq0$. 
In addition it follows in this case $a_1=0$ that 
(ii) $\chi_{q,r}(\a)=Q_qR_r+a_2^3Q_{q-2}R_{r-1}$ for $q,r\geq0$; 
(iii) $Q_q=-a_0a_2Q_{q-2}-a_2^4Q_{q-4}$ for $q\geq1$;
and (iv) $R_r=a_0R_{r-1}-a_0a_2^2R_{r-3}+a_2^4R_{r-4}$ for $r\geq1$. 

As a further special case, we have $\psi_{0,r}(\a)=\chi_{0,r}(\a)-a_2^2\chi_{0,r-2}=R_r-a_2^2R_{r-2}$ from which follow
the first part of (\ref{eqn-oon-psi-AC}), namely the fact that $\psi_{0,0}(\a)=1$ and $\psi_{0,1}=a_0$,
as well as the recurrence relation $\psi_{0,r}(\a)=a_0\psi_{0,r-1}(\a)-a_2^2\psi_{0,r-2}(\a)$ which is a direct consequence of (iv) for $r\geq3$
and also applies in the case $r=2$ for which $\psi_{0,2}=a_0^2-a_2^2$.

If $q=2k+1$ with $k\geq0$, the required result $\psi_{2k+1,r}(\a)=(-a_2)^{r+k+1}\psi_{0,k}(\a)$ may be established as follows.
It should first be noted that from (iii) that $Q_q=0$ not only for $q=1$ but for all $q=2k+1$ with $k\geq0$.
Then (i) and (ii) imply that $\psi_{2k+1,r}(\a)=(-a_2)^{r+1}Q_{2k}$ for all $k\geq0$.
so it simply has to be shown the $Q_{2k}=(-a_2)^k\psi_{0,k}(\a)$ for all $k\geq0$. 
That it is true for $k=0$, $1$ and $2$ as can be seen from the first part of (\ref{eqn-oon-psi-AC}) together with the fact that 
$Q_0=1$, $Q_2=-a_0a_2$ and $Q_4=(a_0^2-a_2^2)a_2^2$. Moreover, (\ref{eqn-oon-psi-AC}) would imply that
\begin{equation}\label{eqn-psi0k}
  (-a_2)^k\psi_{0,k}(\a)=-a_0a_2(-a_2)^{k-1}\psi_{0,k-1}(\a)-a_2^4(-a_2)^{k-2}\psi_{0,k-2}(\a)~~\mbox{for $k\geq3$}\,.
\end{equation}
This recurrence relation for $(-a_2)^k\psi_{0,k}(\a)$ coincides with that of $Q_{2k}$ as given by (iii), so that
by induction on $k$ the required identity between $(-a_2)^k\psi_{0,k}(\a)$ and $Q_{2k}$ is valid for all $k\geq0$.
  
Now we consider the case $q=2k$ with $k\geq0$ for which the required result is $\psi_{2k,r}(\a)=(-a_2)^k\psi_{0,k+r}(\a)$.
This is trivially true for $k=0$. For $k=1$ we have from (i) and (ii) that
\begin{equation}\label{eqn-psi2r}
  \psi_{2,r}(\a)=-a_0a_2R_r+a_2^3R_{r-1}+a_0a_2^3R_{r-2}-a_2^5R_{r-3}~~\mbox{for $r\geq0$}\,.
\end{equation}
Again from (i) and (ii) we have $\psi_{0,r+1}=R_{r+1}-a_2^2R_{r-1}$, and the recurrence relation (v) applied 
to $R_{r+1}$ then yields
\begin{equation}\label{eqn-psi0r1}
  (-a_2)\psi_{0,r+1}=(-a_2)(a_0R_r-a_0a_2^2R_{r-2}+a_2^4R_{r-3})-(-a_2)a_2^2R_{r-1}~~\mbox{for $r\geq0$}\,,
\end{equation}
but this equals $\psi_{2,r}(\a)$ as given in (\ref{eqn-psi2r}), as required.

This is true in the case $r=0$ since (i) and (ii) give $\psi_{2k,0}(\a)=Q_{2k}$, but as shown above $Q_{2k}=(-a_2)^k\psi_{0,k}(\a)$.
The fact that $Q_{2k-1}=0$ for all $k\geq0$ implies that from (i) and (ii) we have
$\psi_{2k,r}(\a)= Q_{2k}(R_r-a_2^2R_{r-2})+a_2^3Q_{2k-2}(R_{r-1}-a_2^2R_{r-3})$ for all $k\geq0$ and $r\geq1$. 
It then follows from the recurrence relations (iii) for $Q_q$ with $q=2k\geq1$ and $q=2k-2\geq1$ that 
$\psi_{2k,r}(\a)=-a_0a_2\psi_{2k-2,r}(\a)+a_2^4\psi_{2k-4,r}(\a)$ for $k\geq2$ and $r\geq1$.
Given the validity of the cases $k=0$ and $k=1$, this allows us to see by induction that for $k\geq2$ and $r\geq1$
\begin{equation}
\begin{array}{rcl}
\psi_{2k,r}(\a)&=&(-a_0a_2)(-a_2)^{k-1}\psi_{0,k-1+r}(\a)+a_2^4(-a_2)^{k-2}\psi_{0,k-2+r}(\a)\cr\cr
               &=&(-a_2)^k(a_0\psi_{0,k-1+r}(\a)+a_2^2\psi_{0,k-2+r}(\a)=(-a_2)^k\psi_{0,k+r}(\a)\,,\cr

\end{array}
\end{equation}
exactly as required, where the final step follows from (\ref{eqn-psi0k}) since $k+r\geq3$. 
\qed

Some results obtained from Theorem~\ref{The-oon} and Corollary~\ref{Cor-oon-AB-AC} for various specific values 
of $\a=(a_0,1)$, $(a_0,0,1)$, $(a_0,1,-1,1)$ and $(a_0,1,1)$ can be found in 
Tables~\ref{Tab-phi-values}, \ref{Tab-oon-psi-A01}, \ref{Tab-oon-psi-A11}and \ref{Tab-oon-psi-A-11},
respectively.


\section{Even orthogonal group identities}\label{sec-eon}

To every partition $\lambda$ of length $\ell(\lambda)\leq n$ there corresponds an irreducible representation of $O(2n)$
of highest weight $\lambda$. On restriction to $SO(2n)$ such a representation remains irreducible if $\ell(\lambda)<n$ 
but decomposes into a sum of two irreducible representations of $SO(2n)$ of highest weights $\lambda_+=(\lambda_1,\ldots,\lambda_{n-1},\lambda_n)$ and 
$\lambda_-=(\lambda_1,\ldots,\lambda_{n-1},-\lambda_n)$. Accordingly, in terms of characters, we have 
\begin{equation}\label{eqn-O2n-SO2n}
   \ch^{O(2n)}_\lambda(\x,\ov\x)=\begin{cases}  
	             \ch^{SO(2n)}_\lambda(\x,\ov\x)&\mbox{if $\ell(\lambda)<n$};\cr
	  \ch^{SO(2n)}_{\lambda_+}(\x,\ov\x)+\ch^{SO(2n)}_{\lambda_-}(\x,\ov\x)&\mbox{if $\ell(\lambda)=n$},\cr
		  \end{cases}
\end{equation}
where the characters of $O(2n)$ have been evaluated on group elements with eigenvalues $(\x,\ov\x)$
which necessarily belong to the subgroup $SO(2n)$. Bearing this notation in mind, we have 
\begin{Theorem}\label{The-eon}
For all $\a=(a_0,a_1,a_2)$ 
\begin{equation}\label{eqn-eon-gf}
\prod_{i=1}^n \left(\,a_0+a_1(x_i+\ov{x}_i)+a_2(x_i^2+\ov{x}_i^2)\,\right) 
  = \sum_{p=0}^n \sum_{q=0}^{n-p} \delta_{r,n-p-q} \ a_2^p\ \psi_{q,r}(\a)\ \ch^{O(2n)}_{(2^p,1^q,0^r)}(\x,\ov\x)\,,
\end{equation}
where
\begin{equation}\label{eqn-eon-psi}
  \psi_{q,r}(\a)=\ds\sum_{t=0}^r \chi_{q,r-t}(\a)\,S_t + \chi_{q-1,0}(1\!-\!\delta_{r0})a_1(-a_2)^{r}+\chi_{q-2,0}(-\delta_{r,0}\,a_2^2+\delta_{r,1}\,a_2^3)\,, 
\end{equation}
with $\chi_{q,r}$ defined as in (\ref{eqn-spn-chi}) in terms of $Q_q$ and $R_r$ given by (\ref{eqn-spn-Q}) and (\ref{eqn-spn-R}), respectively, 
and with
\begin{equation}\label{eqn-eon-S}
\begin{array}{l}
    S_0=1,~~S_1=0,~~S_2=a_0a_2-a_1^2,~~S_3=-a_0a_2^2+2a_1^2a_2-2a_2^3,~~S_4=-2a_1^2a_2^2+a_2^4\cr\cr
		~~\mbox{and}~~S_t=-2a_1^2(-a_2)^{t-2}~~\mbox{for $t\geq5$}\,.\cr
\end{array}
\end{equation}
\end{Theorem}

\noindent{\bf Proof}:
The product on the left hand side of (\ref{eqn-eon-gf}) is clearly invariant under permutations of $x_i$ for $i=1,2,\ldots,n$
and under sign changes of one or more of the exponents of the $x_i$. This symmetry is larger than that of $W_{SO(2n)}$
which only includes even numbers of changes of the exponents of the $x_i$, and ensures that this product is Weyl 
group invariant. Moreover, its expansion as a sum of terms of the form $\x^\kappa$ with $\kappa_j\in\{\ov2,\ov1,0,1,2\}$ 
for $j=1,2,\ldots,n$ is such that each $\kappa$ lies in $\Lambda_{SO(2n)}$. It follows that Proposition~\ref{Pro-WG-invariance}
applies, so that the expansion can be re-expressed as sum of characters $\ch^{SO(2n)}_\lambda(\x,\ov\x)$ with $\lambda\in\Lambda^+$.

In order to identify the $\lambda$ that can appear and the relevant weighting inherited from that of $a(\kappa)$
one proceeds as in both the symplectic and odd orthogonal case by exploiting the transformations and weightings of Table~\ref{Tab-seq-spn}
that arise from permutations of the components of subsequences of $\kappa$. 

There are two significant differences in the even orthogonal case, first the weights of the form $(2,\ldots,2,\ov2)$ 
and $(2,\ldots,2,1,\ldots,1,\ov1)$ lie in $\Lambda_{SO(2n)}^+$. However, the symmetry of product on the left hand side of (\ref{eqn-eon-gf}) 
with respect to sign changes of the exponents of $x_i$ for all $i=1,2,\ldots,n$ is sufficient to ensure that in its
expansion and the subsequent standardisation of $\kappa$ the multiplicities of the terms with exponents $(2^n)$ and $(2^{n-1},\ov2)$
with $n>0$ must be identical, as must those with exponents $(2^p,1^{n-p})$ and $(2^p,1^{n-p-1},\ov1)$ with $n>p>0$. This is what leads to  
the coefficients on the right hand side of (\ref{eqn-eon-gf}) being the same for $\ch^{SO(2n)}_{(2^n)}(\x,\ov\x)$ and 
$\ch^{SO(2n)}_{(2^{n-1},\ov2)}(\x,\ov\x)$, thereby justifying the fact their contribution has been gathered together
in the terms in $\ch^{O(2n)}_{(2^n)}(\x,\ov\x)$, with a similar argument applying to the origin of terms in 
$\ch^{O(2n)}_{(2^p,1^{n-p}}(\x,\ov\x)$. 

The second difference is that in the even orthogonal case the
dot action of $w_{\alpha_n}$ on $\kappa$ with $\alpha_n=\epsilon_{n-1}+\epsilon_n$ now transforms the last two 
non-standard components of $\kappa$ while leaving all other components unchanged
as shown below:  
\begin{equation}\label{eqn-eon-wn}
\begin{array}{|l|lllll|}
\hline
(\ldots,\kappa_{n-1},\kappa_n)&(\ldots,\ov2,\ov2)&(\ldots,\ov2,0)&(\ldots,\ov1,\ov2)&(\ldots,\ov1,\ov1)&(\ldots,0,\ov2)\cr
(\ldots,\!-\!\kappa_n\!-\!1,\!-\!\kappa_{n-1}\!-\!1)&(\ldots,1,1)&(\ldots,\ov1,1)&(\ldots,1,0)&(\ldots,0,0)&(\ldots,1,\ov1)\cr
\hline
(\ldots,\kappa_{n-1},\kappa_n)&(\ldots,\ov2,1)&(\ldots,\ov1,0)&(\ldots,0,\ov1)&(\ldots,1,\ov2)&(\ldots,\ov2,\ov1)\cr
(\ldots,\!-\!\kappa_n\!-\!1,\!-\!\kappa_{n-1}\!-\!1)&(\ldots,\ov2,1)&(\ldots,\ov1,0)&(\ldots,0,\ov1)&(\ldots,1,\ov2)&(\ldots,0,1)\cr
\hline
\end{array}
\end{equation}
As a result terms with $\kappa=(\ldots,\ov2,1)$, $(\ldots,\ov1,0)$, $(\ldots,0,\ov1)$ or $(\ldots,1,\ov2)$ make no contribution. 
The same is true of $(\ldots,\ov2,\ov1)$ since it reduces to minus that of $(\ldots,0,1)$ which is zero, by virtue of (\ref{eqn-spn-wi}). 
In addition that for $\kappa=(\ldots,\ov2,0)$ is equivalent to minus that of $(\ldots,\ov1,1)$ which is in turn equivalent 
to minus that of $(\ldots,0,0)$ by virtue of~(\ref{eqn-spn-wi}), but with differing weights $a_0a_2$, $-a_1^2$ and $a_0^2$, respectively.

It is easy confirm the identity of coefficients of $\lambda=(2^n)$ and $(2^{n-1},\ov2)$ since there are no 
$\kappa$'s with components restricted to $\{\ov2,\ov1,0,1,2\}$ that can lead under the dot action of the Weyl group to 
$(2^n)$ and $(2^{n-1},\ov2)$ other than $(2^n)$ and $(2^{n-1},\ov2)$ themselves and these carry the same weight, namely $a_2^n$. 
The same argument applies to the case of $\lambda=(2^{n-1},1)$ and $(2^{n-1},\ov1)$ but the situation involving 
$\lambda=(2^p,1^{n-p})$ and $(2^p,1^{n-p-1},\ov1)$, that is $\lambda=(\ldots,1,1)$ and $(\ldots,1,\ov1)$
with $\ldots=(2^p,1^{n-p-2})$ for $0\leq p\leq n-2$ is not quite so obvious. 
However, as can be seen from (\ref{eqn-spn-wi}) and (\ref{eqn-eon-wn}), the contributions from $\kappa=(\ldots,\ov2,\ov2)$ 
and $(\ldots,\ov2,2)$ to $\lambda=(\ldots,1,1)$ and $(\ldots,1,\ov1)$ are identical, each carrying a weight $-a_2^2$. 
Similarly, both $\kappa=(\ldots,0,\ov2)$ and $(\ldots,0,2)$, carry identical weights $-a_0a_2$, 
with both $\kappa=(\ldots,1,1)$ and $(\ldots,1,\ov1)$ carrying weight $a_1^2$, thereby confirming in each case 
the equality of the coefficients of $\lambda=(\ldots,1,1)$ and $(\ldots,1,\ov1)$, as already anticipated above.

Combining these observations about the final or terminating entries in any sequence $\kappa$ with those of (\ref{eqn-spn-wi}) one arrives 
at the list of transformations of terminal subsequences given in Table~\ref{Tab-seq-eon} in which $t$ may be any positive integer. 
\begin{table}[ht]
\begin{center}
\begin{tabular}{|l|l|l|}
\hline
$\sigma$&$\tau=w\cdot\sigma$&$\sgn(w)\,a(\sigma)$ \cr
\hline
\hline
$(\cdots,\ov2,\ov2)$&$(\cdots,1,1)$&$-a_2^2$\cr
\hline
$(\cdots,\ov2,2,\ov2)$&$(\cdots,1,1,0)$&$+a_2^3$\cr
\hline
\hline
$(\cdots,\ov1,\ov2)$&$(\cdots,1,0)$&$-a_1a_2$\cr
\hline
$(\cdots,(\ov2,2)^t,\ov1)$&$(\cdots,1,0^{2t})$&$+a_1a_2^{2t}$\cr
$(\cdots,(\ov2,2)^t,\ov1,\ov2)$&$(\cdots,1,0^{2t+1})$&$-a_1a_2^{2t+1}$\cr
\hline
\hline
$(\cdots,\ov2,0)$&$(\cdots,0,0)$&$+a_0a_2$\cr
$(\cdots,\ov1,\ov1)$&$(\cdots,0,0)$&$-a_1^2$\cr
\hline
$(\cdots,\ov2,0,\ov2)$&$(\cdots,0,0,0)$&$-a_0a_2^2$\cr
$(\cdots,\ov1,\ov1,\ov2)$&$(\cdots,0,0,0)$&$+a_1^2a_2$\cr
$(\cdots,\ov2,\ov1,\ov1)$&$(\cdots,0,0,0)$&$+a_1^2a_2$\cr
$(\cdots,\ov2,\ov2,\ov2)$&$(\cdots,0,0,0)$&$-a_2^3$\cr
$(\cdots,\ov2,\ov2,2)$&$(\cdots,0,0,0)$&$-a_2^3$\cr
\hline
$(\cdots,\ov1,\ov2,2,\ov1)$&$(\cdots,0,0,0,0)$&$-a_1^2a_2^2$\cr
$(\cdots,\ov2,1,\ov1,\ov2)$&$(\cdots,0,0,0,0)$&$-a_1^2a_2^2$\cr
$(\cdots,\ov2,\ov2,2,\ov2)$&$(\cdots,0,0,0,0)$&$+a_2^4$\cr
\hline
$(\cdots,\ov1,(\ov2,2)^t,\ov1,\ov2)$&$(\cdots,0^{2t+3})$&$+a_1^2a_2^{2t+1}$\cr
$(\cdots,\ov2,1,(\ov2,2)^t,\ov1)$&$(\cdots,0^{2t+3})$&$+a_1^2a_2^{2t+1}$\cr
$(\cdots,\ov1,(\ov2,2)^{t+1},\ov1)$&$(\cdots,0^{2t+4)})$&$-a_1^2a_2^{2t+2}$\cr
$(\cdots,\ov2,1,(\ov2,2)^t,\ov1,\ov2)$&$(\cdots,0^{2t+4)})$&$-a_1^2a_2^{2t+4}$\cr
\hline
\hline
\end{tabular}
\end{center}
\medskip
\caption{Elementary terminating subsequences $\tau=w(\sigma)$ of $\lambda$ along with $\sgn(w)$ and the contribution of $a(\sigma)$ to $a(\kappa)$ 
in the case $SO(2n)$.} 
\label{Tab-seq-eon}
\end{table}

As in the case of the symplectic and odd orthogonal groups the transformations of Table~\ref{Tab-seq-spn} 
imply the validity of the recurrence relations for $Q_q$, $R_r$ and $\chi_{q,r}(\a)$,
so those of Table~\ref{Tab-seq-eon} imply the necessity of defining $\psi_{q,r}(\a)$ as in (\ref{eqn-eon-psi}) in order to include 
all those contributions leading to $\lambda$ of the form $(\cdots,0^s)$, $(\cdots,1,0^r)$ and $(\cdots,1,1,0^r)$ with the $\cdots$ signifying all terms
enumerated by $\chi_{q,r-s}(\a)$, $\chi_{q-1,0}(\a)$ and $\chi_{q-2,0}(\a)$, respectively.
\qed

Special cases of Theorem~\ref{The-eon} include the following
\begin{Corollary}\label{Cor-eon-AB-AC}
\begin{align}
\ds \prod_{i=1}^n (a_0+a_1(x_i+\ov{x}_i)) &= \sum_{q=0}^{n}\,\delta_{r,n-q}\,a_1^q\,\phi_r(\a)\,\ch^{O(2n)}_{(1^q,0^r)}(\x,\ov\x)\,;\label{eqn-eon-AB}\\
\ds \prod_{i=1}^n (a_0+a_2(x_i^2+\ov{x}_i^2)) &= \sum_{p=0}^n \sum_{q=0}^{n-p}\,\delta_{r,n-p-q}\,a_2^p\,\psi_{q,r}(\a)\,\ch^{O(2n)}_{(2^p,1^q,0^r)}(\x,\ov\x)\,,\label{eqn-eon-AC}
\end{align}
with
\begin{equation}\label{eqn-eon-phi-AB}
\begin{array}{l}
\!\!\phi_0(\a)=1, \phi_1=a_0, \phi_2=a_0^2-2a_1^2~\mbox{and $\phi_r(\a)=a_0\,\phi_{r-1}(\a)-a_1^2\,\phi_{r-2}(\a)$ for $r\geq3$},\cr
\end{array}
\end{equation}
and
\begin{equation}\label{eqn-eon-psi-AC}
\begin{array}{rl}
\psi_{q,r}(\a)&=0~~\mbox{if $q<0$, $q=1\mod2$ or $r<0$}\,,\cr
\psi_{0,0}(\a)&=1,~~\psi_{0,1}(\a)=a_0,~~\psi_{0,2}(\a)=a_0^2+a_0a_2\,,\cr
\psi_{0,3}(\a)&=a_0^3+a_0^2a_2-2a_0a_2^2-2a_2^3\,,
\psi_{0,4}(\a)=a_0^4+a_0^3a_2-3a_0^2a_2^2-2a_0a_2^3+2a_2^4\,,\cr
\psi_{0,r}(\a)&=a_0\psi_{0,r-1}(\a)-a_0a_2^2\psi_{0,r-3}(\a)+a_2^4\psi_{0,r-4}(\a)~~\mbox{for $r\geq5$}\,,\cr
\psi_{2,0}(\a)&=-a_0a_2-a_2^2,~~\psi_{2,1}(\a)=-a_0^2a_2+2a_2^3\,,\cr
\psi_{2,r}(\a)&=-a_0a_2\psi_{0,r}(\a)+a_2^3\psi_{0,r-1}(\a)~~\mbox{for $r\geq2$}\,,\cr
\psi_{q,r}(\a)&=-a_0a_2\psi_{q-2,r}(\a)-a_2^4\psi_{q-4,r}(\a)~~\mbox{for $q=0\mod2$ and $q\geq4$ and $r\geq0$}\,.
\end{array}
\end{equation}
\end{Corollary}

\noindent{\bf Proof}:
The proofs of both (\ref{eqn-eon-phi-AB}) and (\ref{eqn-eon-psi-AC}) proceed in the same way as those of (\ref{eqn-spn-phi-AB}) and (\ref{eqn-spn-psi-AC}).
The only differences are those applying to low values of $q$ and $r$ which have been established by explicit calculation.
\qed

Some results obtained from Theorem~\ref{The-eon} and Corollary~\ref{Cor-eon-AB-AC} for various specific values 
of $\a=(a_0,1)$, $(a_0,0,1)$ and $(a_0,1,1)$ are shown in 
Tables~\ref{Tab-phi-values}, \ref{Tab-eon-psi-A01} and \ref{Tab-eon-psi-A11},
respectively.

\section{Expansions in terms of spin characters} \label{sec-spin}

In the case of $SO(2n+1)$ and $SO(2n)$ consideration has been given so far to cases for which
the products only involve weight space vectors whose components are all integers, even though
the sets of all finite dimensional irreducible representations of these groups also includes those whose weights involve
vectors whose components are all half odd integers. Such representations are often referred to as spin representations
since they are faithful irreducible representations of the spin covering groups of $SO(2n+1)$ and $SO(2n)$.
For every partition $\mu$ of length $\ell(\mu)\leq n$ there exists an irreducible representation of
of $SO(2n+1)$ of highest weight $\lambda=\Delta+\mu$, where $\Delta=(\1,\1,\ldots,\1)$ has length $n$.
The same is true of $O(2n)$, but on restriction to $SO(2n)$ this representation decomposes into a sum 
of two irreducible representations of highest weights $\lambda_+=(\lambda_1,\lambda_2,\ldots,\lambda_{n-1},\lambda_n)$
and $\lambda_-=(\lambda_1,\lambda_2,\ldots,\lambda_{n-1},-\lambda_n)$

The characters of $SO(2n+1)$ and of $O(2n)$ of highest weights $\lambda=\Delta+\mu$ are denoted,
when restricted to group elements with eigenvalues $(\x,\ov\x,1)$ and $(\x,\ov\x)$, by
$\ch_{\Delta+\mu}^{SO(2n+1)}(\x,\ov\x,1)$ and $\ch_{\Delta+\mu}^{O(2n)}(\x,\ov\x)$, respectively. 
In particular the basic spin representations of both $SO(2n+1)$ and $O(2n)$ have highest weight 
$\Delta=(\1,\1,\ldots,\1)$ and identical characters
\begin{equation}\label{eqn-Delta}
   \ch_\Delta^{SO(2n+1)}(\x,\ov\x,1)=\ch_\Delta^{O(2n)}(\x,\ov\x)=\prod_{i=1}^n (x_i^\1+x_i^{-\1})\,.
\end{equation}

The omission of identities involving spin characters can be remedied through a consideration of products of the
form: 
\begin{equation}\label{eqn-oon-spin-axkappa}
   \prod_{i=1}^n \sum_{k=1}^m\, a_{k+\1} (x_i^{k+\1}+x_i^{-k-\1}) = \sum_{\kappa}\, a(\kappa)\, \x^{\kappa} 
\end{equation}
where in the expansion on the right the components of $\kappa$ are now all half-integral, or more precisely half an odd integer. 

Proposition~\ref{Pro-WG-invariance} still applies for both $SO(2n+1)$ and $SO(2n)$, as does the dot action of
the generators of the corresponding Weyl groups on $\kappa$ that is set out in Table~\ref{Tab-wkappa}. 
By way of an example, if one restricts oneself to the case $m=1$ then the expansion of the above product 
in terms of characters of $SO(2n+1)$ is given in

\begin{Theorem}\label{The-oon-spin-gf-AB}~~For $\a=(a_\1,a_\3)$
\begin{equation}\label{eqn-oon-spin-gf-AB}
    \prod_{i=1}^n \left(\,a_{\1}(x_i^{\1}+x_i^{-\1})+a_{\3}(x_i^{\3}+x_i^{-\3})\,\right) 
		= \sum_{q=0}^n\, \delta_{r,n-q}\, a_{\3}^q\, \phi_r(\a)\, \ch^{SO(2n+1)}_{(\3^q,\1^r)}(\x,\ov\x,1)\,,
\end{equation}
where 
\begin{equation}\label{eqn-oon-spin-phi-AB}
\begin{array}{l}
\phi_0(\a)=1, \phi_1(\a)=a_{\1}-a_{\3}, \phi_2(\a)=a_{\1}^2-2a_{\1} a_{\3}~~\mbox{and}~~ \cr
\phi_r(\a)=a_{\1}\phi_{r-1}(\a)-a_{\1} a_{\3}\phi_{r-2}(\a)+a_{\3}^3\phi_{r-3}(\a)~~\mbox{if $r\geq3$}.\cr
\end{array}
\end{equation}
\end{Theorem}

\noindent{\bf Proof}: 
The half odd integer versions of (\ref{eqn-spn-wi}) and (\ref{eqn-oon-wn})
appropriate to $SO(2n+1)$ follow from the formulae in Table~\ref{Tab-wkappa}. 
They take the form 
\begin{equation}\label{eqn-oon-spin-wi-AB}
\begin{array}{|l|llllll|}
\hline
(\kappa_i,\kappa_{i+1})&(\ov\3,\ov\1)&(\ov\3,{\1})&(\ov\3,{\3})&(\ov\1,{\1})&(\ov\1,{\3})&({\1},{\3})\cr
(\kappa_{i+1}-1,\kappa_{i}+1)&(\ov\3,\ov\1)&(\ov\1,\ov\1)&({\1},\ov\1)&(\ov\1,{\1})&({\1},{\1})&({\1},{\3})\cr
\hline
\end{array}
\end{equation}
and 
\begin{equation}\label{eqn-oon-spin-wn-AB}
\begin{array}{|l|ll|}
\hline
(\ldots,\kappa_n)&(\ldots,\ov{{\3}})&(\ldots,\ov{{\1}})\cr
(\ldots,-\kappa_n-1)&(\ldots,{\1})&(\ldots,\ov{{\1}})\cr
\hline
\end{array}
\end{equation} 
These lead directly to the following transformations of elementary and terminating subsequences: 
\begin{equation}\label{eqn-oon-spin-sigma-tau-wi-AB}
\begin{array}{|l|l|l|}
\hline
\sigma&\tau=w\cdot\sigma&\sgn(w)\,a(\sigma) \cr
\hline\hline(\ldots,\1,\ldots)&(\ldots,\1,\ldots)&+a_\1\cr
(\ldots,\ov\1,\3,\ldots)&(\ldots,\1,\1,\ldots)&-a_\1 a_\3 \cr
(\ldots,\ov\3,\3,\3,\ldots)&(\ldots,\1,\1,\1,\ldots)&+a_\3^3\cr
\hline
\end{array}
\end{equation}
and
\begin{equation}\label{eqn-oon-spin-sigma-tau-wn-AB}
\begin{array}{|l|l|l|}
\hline
\sigma&\tau=w\cdot\sigma&\sgn(w)\,a(\sigma) \cr
\hline
\hline
(\ldots,\ov\3)&(\ldots,\1)&-a_\3\cr
\hline
\end{array}
\end{equation}
These then imply the validity of the required recurrence relations (\ref{eqn-oon-spin-phi-AB}).
\qed

Similarly, in the case of $SO(2n)$ the comparable result takes the form

\begin{Theorem}\label{Tab-eon-spin-AB}~~For $\a=(a_\1,a_\3)$
\begin{equation}\label{eqn-eon-spin-gf-AB}
    \prod_{i=1}^n \left(\,a_{\1}(x_i^{\1}+x_i^{-\1})+a_{\3}(x_i^{\3}+x^{-\3})\,\right) 
		= \sum_{q=0}^n\, \delta_{r,n-q}\, a_{\3}^q\, \phi_r(\a)\, \ch^{O(2n)}_{(\3^q,\1^r)}(\x,\ov\x)\,,
\end{equation}
where 
\begin{equation}\label{eqn-eon-spin-phi-AB}
\begin{array}{l}
\phi_0(\a)=1, \phi_1(\a)=a_{\1}, \phi_2(\a)=a_{\1}^2-a_{\1} a_{\3}-a_{\3}^2~~\mbox{and}~~\cr\cr
\phi_r(\a)=a_{\1}\phi_{r-1}(\a)-a_{\1} a_{\3}\phi_{r-2}(\a)+a_{\3}^3\phi_{r-3}(\a)~~\mbox{if $r\geq3$}\cr
\end{array}
\end{equation}
\end{Theorem}

\noindent{\bf Proof}:
The $SO(2n)$ formulae in Table~\ref{Tab-wkappa} once again yield the transformations (\ref{eqn-oon-spin-wi-AB})
but those are now to be augmented by 
\begin{equation}\label{eqn-eon-spin-wn-AB}
\begin{array}{|l|llll|}
\hline
(\ldots,\kappa_{n-1},\kappa_n)&(\ldots,\ov{{\3}},\ov{{\3}})&(\ldots,\ov{{\1}},\ov{{\3}})&(\ldots,{\1},\ov{{\3}})&(\ldots,\ov{{\1}},\ov{{\1}})\cr
(\ldots,-\kappa_n-1,-\kappa_{n-1}-1)&(\ldots,{\1},{\1})&(\ldots,{\1},\ov{{\1}})&(\ldots,{\1},\ov{{\3}})&(\ldots,\ov{{\1}},\ov{{\1}})\cr
\hline
\end{array}
\end{equation} 
These lead to the transformations of elementary subsequences given in (\ref{eqn-oon-spin-sigma-tau-wi-AB})
and those of terminating subsequences given by
\begin{equation}\label{eqn-eon-spin-sigma-tau-wn-AB}
\begin{array}{|l|l|l|}
\hline
\sigma&\tau=w\cdot\sigma&\sgn(w)\,a(\sigma)\cr
\hline\hline
(\ldots,\ov{\3},\ov{\3})&(\ldots,\1,\1)&-a_\3^2\cr
(\ldots,\ov{\1},\ov{\3})&(\ldots,\1,\ov{\1})&-a_\1 a_\3\cr
\hline
\end{array}
\end{equation}
Taken together these lead in turn to the required recurrence relations (\ref{eqn-eon-spin-phi-AB}).
\qed

Although the $SO(2n+1)$ spin character $m=1$ recurrence relation (\ref{eqn-oon-spin-phi-AB}) 
appears to be somewhat different from the corresponding $Sp(2n)$ recurrence relation~(\ref{eqn-spn-phi-AB})
it is not difficult to see by explicit calculation that $\phi_r(a_\1,1)$ in the odd orthogonal spin case
is identical with $\phi_r(a_\1-1,1)$ in the symplectic case. On the other hand the $O(2n)$ spin character
$m=1$ recurrence relation (\ref{eqn-eon-spin-phi-AB}) differs from that of the corresponding $SO(2n+1)$ 
recurrence relation~(\ref{eqn-oon-phi-AB}). However, explicit calculation reveals that fact that $\phi_r(a_\1,1)$ 
in the even orthogonal spin case is identical to $(-1)^r\,\phi_r(-a_\1+1,1)$ in the odd orthogonal case.

These coincidences are not, of course, accidental, and will be explained in the next section 
on dual pairs of groups. Moreover, similar coincidences apply for all values of $m$, allowing
all expansion coefficients in spin orthogonal cases to be evaluated, at least in principle, from those
encountered in the symplectic and ordinary, non-spin, odd orthogonal cases. For this reason
the derivation of $m=2$ orthogonal spin character identities by means of recurrence relations
is relegated to Appendix~\ref{Appendix} where it will be seen that a further degree of complexity 
is encountered even in the simplest type of $m=2$ case in which $\a=(a_\1,0,a_\5)$.


\section{Dual pair approach}~\label{sec-dual-pair}

An alternative approach to generating function for series of characters, both ordinary and spin, is to exploit identities ~\cite{Bau,Mor,JM,Has}
that are analogous to the dual Cauchy identity (\ref{eqn-gln-glm}).

\begin{equation}\label{eqn-dual-ids}
\begin{array}{rcl}
\ch_\Delta^{O(4nm)}((\x,\ov\x)\!\!\times\!\!(\y,\ov\y))\!\!
&=&\ds\!\! \sum_{\lambda\in(m^n)}\! \ch_\lambda^{Sp(2n)}(\x,\ov\x)\, \ch_{\tilde\lambda}^{Sp(2m)}(\y,\ov\y); \\
\ch_\Delta^{O(4nm\!+\!2m)}((\x,\ov\x,1)\!\!\times\!\!(\y,\ov\y))\!\!
&=&\ds\!\! \sum_{\lambda\in(m^n)}\! \ch_{\lambda}^{SO(2n\!+\!1)}(\x,\ov\x,1)\, \ch_{\Delta\!+\!\tilde\lambda}^{O(2m)}(\y,\ov\y);\\
\ch_\Delta^{O(4nm)}((\x,\ov\x)\!\!\times\!\!(\y,\ov\y))\!\!
&=&\ds\!\! \sum_{\lambda\in(m^n)}\! \ch_\lambda^{O(2n)}(\x,\ov\x)\, \ch_{\tilde\lambda}^{O(2m)}(\y,\ov\y);\\
\ch_\Delta^{O(4nm\!+\!2n)}((\x,\ov\x)\!\!\times\!\!(\y,\ov\y,1)))\!\!
&=&\ds\!\! \sum_{\lambda\in(m^n)}\!\! \ch_{\Delta\!+\!\lambda}^{O(2n)}(\x,\ov\x)\, \ch_{\tilde\lambda}^{SO(2m\!+\!1)}(\y,\ov\y,1);\\
\ch_\Delta^{O(4nm\!+\!2n\!+\!2m\!+\!1)}((\x,\ov\x,1)\!\!\times\!\!(\y,\ov\y,1))\!\!
&=&\ds\!\! \sum_{\lambda\in(m^n)}\!\! \ch_{\Delta\!+\!\lambda}^{SO(2n\!+\!1)}(\x,\ov\x,1)\, \ch_{\Delta\!+\!\tilde\lambda}^{SO(2m\!+\!1)}(\y,\ov\y,1),\\
\end{array}
\end{equation}
where $(\x,\ov\x)\!\times\!(\y,\ov\y)=(\ldots,x_iy_j,x_i\ov{y}_j,\ov{x}_iy_j,\ov{x}_i\ov{y}_j,\ldots)$,
while $(\x,\ov\x,1)\!\times\!(\y,\ov\y)$ and $(\x,\ov\x)\!\times\!(\y,\ov\y,1)$ include additional terms
$(\ldots,y_j,\ov{y}_j,\ldots)$ and $(\ldots,x_i,\ov{x}_i,\ldots)$, respectively,
and $(\x,\ov\x,1)\!\times\!(\y,\ov\y,1)$ includes both of these together with $1$, all with $i=1,2,\ldots,n$
and $j=1,2,\ldots,m$.
These all owe their origin to the duality of branching rules for the basic spin representation
of an orthogonal group, $O(NM)$, to subgroups consisting of a direct product of either a pair of two symplectic subgroups, $Sp(N)\times Sp(M)$,
or a pair of two orthogonal subgroups, $O(N)\times O(M)$, with the action of the constituent subgroups in each pair mutually commuting 
and centralising one another in the basic spin representation. 

Thanks to (\ref{eqn-Delta}) the left hand sides of these identities (\ref{eqn-dual-ids}) may be expressed as 
products. They then yield the following expansions
\begin{align}
&\ds \!\prod_{i=1}^n \prod_{j=1}^m(x_i\!+\!\ovx_i\!+\!y_j\!+\!\ovy_j)
\!=\!	\ds \!\!\sum_{\lambda\in(m^n)}\!\! \ch_\lambda^{Sp(2n)}(\x,\ov\x)\, \ch_{\tilde\lambda}^{Sp(2m)}(\y,\ov\y)\,;\label{eqn-gf-spn-spm}\\
&\ds \!\prod_{i=1}^n \prod_{j=1}^m(x_i\!+\!\ovx_i\!+\!y_j\!+\!\ovy_j)
\!=\!	\ds \!\!\sum_{\lambda\in(m^n)}\!\! \ch_\lambda^{SO(2n\!+\!1)}(\x,\ov\x,1)(\!-1)^{|\lambda|}\,\ch_{\tilde\lambda}^{SO(2m\!+\!1)}(\!-\y,\!-\ov\y,1)\,;\label{eqn-gf-oon-oom}\\	
&\ds \!\prod_{i=1}^n \prod_{j=1}^m(x_i\!+\!\ovx_i\!+\!y_j\!+\!\ovy_j)
=	\ds \!\!\sum_{\lambda\in(m^n)}\!\! \ch_\lambda^{O(2n)}(\x,\ov\x)\, \ch_{\tilde\lambda}^{O(2m)}(\y,\ov\y)\,;\label{eqn-gf-eon-eom}\\
&\ds \!\prod_{i=1}^n \prod_{j=1}^m(x_i\!+\!\ovx_i\!+\!y_j\!+\!\ovy_j)\prod_{i=1}^n (x_i^{\1}\!+\!\ovx_i^{\1})
=	\!\!\sum_{\lambda\in(m^n)}\!\! \ch_{\Delta\!+\!\lambda}^{O(2n)}(\x,\ov\x)\, \ch_{\tilde\lambda}^{SO(2m\!+\!1)}(\y,\ov\y,1)\,;\label{eqn-gf-eon-oom}\\					
&\ds \!\prod_{i=1}^n \prod_{j=1}^m(x_i\!+\!\ovx_i\!+\!y_j\!+\!\ovy_j)\prod_{i=1}^n (x_i^{\1}\!+\!\ovx_i^{\1})
=	\!\!\sum_{\lambda\in(m^n)}\!\! \ch_{\Delta\!+\!\lambda}^{SO(2n\!+\!1)}(\x,\ov\x,1)\, \ch_{\tilde\lambda}^{Sp(2m)}(\y,\ov\y),\label{eqn-gf-oon-spm}	      
\end{align}
where use has been made of the identities 
\begin{equation}\label{eqn-oon-spin-spn}
\begin{array}{l}
\ch_{\Delta\!+\!\tilde\lambda}^{O(2m)}(\y,\ov\y) \big/ \ch_{\Delta}^{O(2m)}(\y,\ov\y) \!=\!\ch_{\tilde\lambda}^{O(2m+1))}(\y,\ov\y,\!-1) 
    \!=\! (-1)^{|\tilde\lambda|}\,\ch_{\tilde\lambda}^{SO(2m+1))}(\!-\y,\!-\ov\y,1);\cr\cr
\ch_{\Delta\!+\!\tilde\lambda}^{SO(2m+1)}(\y,\ov\y,1) \big/ \ch_{\Delta}^{SO(2m+1)}(\y,\ov\y,1)=\ch_{\tilde\lambda}^{Sp(2m)}(\y,\ov\y)\,,
\end{array}
\end{equation}
In the above expansions the left hand sides are Weyl group symmetric products of the type discussed throughout 
this paper. The most striking implication of (\ref{eqn-gf-spn-spm})-(\ref{eqn-gf-oon-spm}) is that, viewed as generating functions for series of characters
of symplectic and orthogonal groups of rank $n$, the expansion coefficients are themselves characters of symplectic and orthogonal groups of rank 
$m$. Not only that, the expansion coefficients in (\ref{eqn-gf-oon-spm}) coincide with those of (\ref{eqn-gf-spn-spm}), while those of 
(\ref{eqn-gf-eon-oom}) can be quite readily calculated from those (\ref{eqn-gf-oon-oom}). 

To make the connection with our previous notation it might be noted that in the cases $m=1$ and $m=2$ 
these products are of the form
\begin{equation}\label{eqn-yj-ak}
\begin{array}{l}
   \ds \prod_{i=1}^n (a_0\!+\!a_1(x_i\!+\!\ovx_i))~~\mbox{with $a_0=y_1\!+\!\ovy_1$ and $a_1=1$};\cr 
	 \ds \prod_{i=1}^n (a_\1(x_i^{\1}\!+\!\ovx_i^{\1})\!+\!a_\3(x_i^{\3}\!+\!\ovx_i^{\3}))~~\mbox{with $a_\1=y_1\!+\!\ovy_1+1$ and $a_\3=1$};\cr
   \ds \prod_{i=1}^n (a_0\!+\!a_1(x_i\!+\!\ovx_i)\!+\!a_2(x_i^2\!+\!\ovx_i^2))\cr
		~~~~~~~~~~\mbox{with $a_0=(y_1\!+\!\ovy_1)(y_2\!+\!\ovy_2)\!+\!2$, $a_1=(y_1\!+\!\ovy_1\!+\!y_2\!+\!\ovy_2)$ and $a_2=1$};\cr
	 \ds \prod_{i=1}^n (a_\1(x_i^{\1}\!+\!\ovx_i^{\1})\!+\!a_\3(x_i^{\3}\!+\!\ovx_i^{\3})\!+\!a_\5(x_i^{\5}\!+\!\ovx_i^{\5}))\cr
		~~~~~~~~~~\mbox{with $a_\1=(y_1\!+\!\ovy_1)(y_2\!+\!\ovy_2)\!+(\!y_1\!+\!\ovy_1\!+\!y_2\!+\!\ovy_2)\!+\!2$}\cr
		~~~~~~~~~~\mbox{\phantom{with}} \mbox{$a_\3=(y_1\!+\!\ovy_1\!+\!y_2\!+\!\ovy_2)\!+\!1$ and $a_\5=1$},
\end{array}
\end{equation}
respectively.

In the same two cases $m=1$ and $m=2$ the right hand sides of (\ref{eqn-gf-spn-spm})-(\ref{eqn-gf-oon-spm})  
are expansions in terms of group characters of highest weight $\lambda$ or $\Delta\!+\!\lambda$ with 
$\lambda=(1^q,0^r)$ and $(2^p,1^q,0^r)$ whose coefficients that are themselves dual group characters
of highest weight $\tilde\lambda$ with $\tilde\lambda=(r)$ and $(q\!+\!r,r)$, respectively.
Provided that $a_1=1$, $a_\3=1$, $a_2=1$ or $a_\5=1$, as appropriate, this enables one to identify the expansion 
parameters $\phi_{r}(\a)$ and $\psi_{q,r}(\a)$ that appear in the Theorems and Corollaries of Sections~\ref{sec-spn}-\ref{sec-spin} 
with the rather simple characters of $Sp(2)$, $SO(3)$ or $O(2)$ if $m=1$, and of $Sp(4)$, $SO(5)$ or $O(4)$ if $m=2$,
This identification does not lend itself to deriving the recurrence relations provided in these Sections, but does enable 
results obtained through the use of the recurrence relations to be checked for various $\a$ for a considerable range of 
values of $p$, $q$ and $r$.  

The coincidence of the expansion coefficients in (\ref{eqn-gf-oon-spm}) and (\ref{eqn-gf-spn-spm})
the close relationship between expansion coefficients in (\ref{eqn-gf-eon-oom}) and (\ref{eqn-gf-oon-oom})
is such that the use of (\ref{eqn-yj-ak}) implies 
\begin{equation}\label{eqn-phi-psi-oon-spn}
\begin{array}{rcl}
  \phi_r^{SO(2n\!+\!1)}(a_\1,1)&=&\phi_r^{Sp(2n)}(a_\1-1,1)\,;\cr\cr
	\phi_r^{O(2n)}(a_\1,1)&=&(-1)^r\,\phi_r^{SO(2n\!+\!1)}(-a_\1+1,1)\,;\cr\cr
  \psi_{q,r}^{SO(2n\!+\!1)}(a_\1,a_\3,1)&=&\psi_{q,r}^{Sp(2n)}(a_\1-a_\3+1,a_\3-1,1)\,;\cr\cr
	\psi_{q,r}^{O(2n)}(a_\1,a_\3,1)&=&(-1)^q\,\psi_{q,r}^{SO(2n\!+\!1)}(a_\1-a_\3+1,-a_\3+1,1)\,.
\end{array}
\end{equation}
Here we have taken the liberty of augmenting $\phi_r(\a)$ and $\psi_{q,r}(\a)$ with superscripts 
indicating the $n$-dependent group in terms of whose characters each expansion takes place.  

The first two identities explain the coincidences in the $m=1$ evaluations of $\phi_r(\a)$ pointed out at
end of Section~\ref{sec-spin} and the second two imply moreover that the complexities of the recurrence relations 
in the $m=2$ evaluation of $\psi_{q,r}(\a)$ may be avoided by using the much simpler recurrence relations of 
Sections~\ref{sec-spn} and~\ref{sec-oon}. This complexity shows itself even the simplest $m=2$ cases
with $\a=(a_1,0,1)$, as is made explicit in Appendix~\ref{Appendix}.

\section{Symplectic and orthogonal group identities}~\label{sec-spn-oon-eon-values}

For the groups $G=Sp(2n)$, $SO(2n+1)$ and $O(2n)$ the $m=1$, $\a=(a_0,a_1)$ expansions
\begin{equation}
  \prod_{i=1}^n (a_0+a_1(x_i+\ovx_i))=\sum_{q=1}^n \delta_{r,n-q}\, a_1^q\,\phi_r(\a) \, \ch_{(1^q,0^r)}^G(\x,\ov\x)
\end{equation}
involve coefficients $\phi_r(\a)$ 
specified recursively in (\ref{eqn-spn-phi-AB}), (\ref{eqn-oon-phi-AB}) and (\ref{eqn-eon-phi-AB}). 
The results can be summarised as follows:
\begin{equation}
\begin{array}{|l|l|l|l|l|}
\hline
G&\phi_0(\a)&\phi_1(\a)&\phi_2(\a)&\phi_r(\a) \mbox{for all $r\geq3$}\cr
\hline
\hline
 Sp(2n)&1&a_0&a_0^2-a_1^2&a_0\,\phi_{r-1}(\a)-a_1^2\phi_{r-2}(\a)\\
\hline
 SO(2n\!+\!1)&1&a_0-a_1&a_0^2-a_0a_1-a_1^2&a_0\,\phi_{r-1}(\a)-a_1^2\phi_{r-2}(\a)\\
\hline
 O(2n)&1&a_0&a_0^2-2a_1^2&a_0\,\phi_{r-1}(\a)-a_1^2\phi_{r-2}(\a)\\
\hline
\end{array} 
\end{equation}
or equivalently,
\begin{equation}
\begin{array}{llcl}
 Sp(2n)&\phi_{r}(\a)&=&[t^r]\, 1/(1-a_0\,t+a_1^2\,t^2)\\ 
 SO(2n\!+\!1)&\phi_{r}(\a)&=&[t^r]\, (1-a_1\,t)/(1-a_0\,t+a_1^2\,t^2)\\ 
 O(2n)&\phi_{r}(\a)&=&[t^r]\, (1-a_1^2\,t^2)/(1-a_0\,t+a_1^2\,t^2)
\end{array} 
\end{equation}
 
More explicitly, for all $\a=(a_0,1)$ we have $\phi_0(\a)=1$ with $\phi_r(\a)$ given for all $r>0$ 
and various values of $a_0$ in Table~\ref{Tab-phi-values}.


\begin{table}[ht]
\begin{center}
\begin{tabular}{|c||c||c||c|c|}
\hline
    &$Sp(2n)$&$SO(2n+1)$&$O(2n)$&\cr
$\a$&$\phi_r(\a)$&$\phi_r(\a)$&$\phi_r(\a)$&\cr
\hline
\hline
$(0,1)$&$\begin{array}{rl} 1&\mbox{if $r=0$}\cr 0&\mbox{if $r=1,3$}\cr -1&\mbox{if $r=2$}\cr \end{array}$
       &$\begin{array}{rl} 1&\mbox{if $r=0,3$}\cr -1&\mbox{if $r=1,2$}\cr \end{array}$
       &$\begin{array}{rl} 2&\mbox{if $r=0$}\cr 0&\mbox{if $r=1,3$}\cr -2&\mbox{if $r=2$}\cr\end{array}$
			 &$ \mod4$\cr
\hline
$(1,1)$&$\begin{array}{rl} 1&\mbox{if $r=0,1$}\cr 0&\mbox{if $r=2,5$}\cr -1&\mbox{if $r=3,4$}\cr \end{array}$
       &$\begin{array}{rl} 1&\mbox{if $r=0,5$}\cr 0&\mbox{if $r=1,4$}\cr -1&\mbox{if $r=2,3$}\cr \end{array}$
       &$\begin{array}{rl} 2&\mbox{if $r=0$}\cr 1&\mbox{if $r=1,5$}\cr -1&\mbox{if $r=2,4$}\cr -2&\mbox{if $r=3$}\cr\end{array}$
			 &$\mod6$\cr
\hline
$(-1,1)$&$\begin{array}{rl} 1&\mbox{if $r=0$}\cr 0&\mbox{if $r=1$}\cr -1&\mbox{if $r=2$}\cr \end{array}$
       &$\begin{array}{rl} 1&\mbox{if $r=0,2$}\cr -2&\mbox{if $r=1$}\cr  \end{array}$
       &$\begin{array}{rl} 2&\mbox{if $r=0$}\cr -1&\mbox{if $r=1,2$}\cr \end{array}$
			 &$\mod3$\cr
\hline
$(2,1)$&$r+1$&$1$&$2$&$\mbox{all $r$}$\cr
\hline
$(-2,1)$&$(-1)^r(r+1)$&$(-1)^r(2r+1)$&$(-1)^r 2$&$\mbox{all $r$}$\cr
\hline
$(3,1)$&$F_{2r+2}$&$F_{2r+1}$&$F_{2r+2}-F_{2r-2}$&$\mbox{all $r$}$\cr
\hline
$(-3,1)$&$(-1)^r F_{2r+2}$&$(-1)^r(F_{2r+2}+F_{2r})$&$(-1)^r(F_{2r+2}-F_{2r-2})$&$\mbox{all $r$}$\cr
\hline
\end{tabular}
\end{center}
\medskip
\caption{The coefficients $\phi_r(\a)$ for various $\a=(a_0,1)$ where $F_k$ is the $k$th Fibonacci number} 
\label{Tab-phi-values}																																								
\end{table}


For $Sp(2n)$ in the special $m=2$ case with $a_1=0$, that is $\a=(a_0,0,a_2)$, the expansion
\begin{equation}
  \prod_{i=1}^n (a_0+a_2(x_i^2+\ovx_i^2))=\sum_{p=1}^n\, \sum_{q=1}^{n-p}\, \delta_{r,n-p-q}\, a_2^p\,\psi_{q,r}(\a) \, \ch_{(2^p,1^q,0^r)}^{Sp(2n}(\x,\ov\x)
\end{equation}
involves coefficients $\psi_{q,r}(\a)$ which may be evaluated recursively
by means of (\ref{eqn-spn-psi-AC}). Values obtained in this way for various $\a=(a_0,0,1)$ are displayed in Table~\ref{Tab-spn-psi-A01}.
These include results obtained previously by Lee and Oh~\cite{LO1}. 

By way of a further $Sp(2n)$ example, for $\a=(a_0,1,1)$ explicit results for the coefficients $\psi_{q,r}(\a)$ appearing in
\begin{equation}\label{eqn-spn-psi-A11}
   \prod_{i=1}^n (a_0+(x_i+\ov{x}_i)+(x_i^2+\ov{x}_i^2))
	     =\sum_{p=0}^{n} \sum_{q=0}^{n-p}\ \delta_{r,n-p-q}\ \psi_{q,r}(\a)\ \ch^{Sp(2n)}_{(2^p,1^q,0^r)}(\x,\ov\x)\,,
\end{equation}
may be obtained from Theorem~\ref{The-spn}. Some of these are displayed in Table~\ref{Tab-spn-psi-A11}. 
Mapping $(a_0,a_1,a_2)$ to $(a_0,-a_1,a_2)$ is equivalent to mapping
$\x$ to $-\x$ and $\ov\x$ to $-\ov\x$. However, in contrast to the $GL(n)$ case, the terms appearing 
in $\ch^{Sp(2n)}_\lambda(\x,\ov\x)$ are not homogeneous in the components of $\x$. Nonetheless for each constituent multinomial in the
components of $\x$ and $\ov\x$ the sum of all the exponents of $x_i$ for $i=1,2,\ldots,n$ is even or odd 
according as $|\lambda|$ is even or odd. This implies that $\psi_{q,r}(a_0,-a_1,1)=(-1)^q\,\psi_{q,r}(a_0,a_1,1)$.
Thus the results of Table~\ref{Tab-spn-psi-A11} may be extended by changing $a_1=1$ to $a_1=-1$ and including an additional 
factor of $(-1)^q$.

For $SO(2n\!+\!1)$ in the case $m=2$ and $\a=(a_0,0,1)$ the coefficients $\psi_{q,r}(\a)$ appearing in the expansion 
\begin{equation}\label{eqn-oon-psi-A01}
   \prod_{i=1}^n (a_0+(x_i^2+\ov{x}_i^2))
	     =\sum_{p=0}^n \sum_{q=0}^{n-p} \delta_{r,n-p q}\ \psi_{q,r}(\a)\ \ch^{SO(2n+1)}_{(2^p,1^q,0^r)}(\x,\ov\x,1)\,,
\end{equation}
may be evaluated from the second part of Corollary~\ref{Cor-oon-AB-AC}. Examples of this type are displayed in Table~\ref{Tab-oon-psi-A01}.

Similarly, in the case $\a=(a_0,1,1)$ the $SO(2n+1)$ expansion coefficients $\psi_{q,r}(\a)$ appearing in 
\begin{equation}\label{eqn-oon-psi-A11}        
\prod_{i=1}^n (a_0+x_i+\ov{x}_i+x_i^2+\ov{x}_i^2) = \sum_{p=0}^n \sum_{q=0}^{n-p} \delta_{r,n-p-q}\ \psi_{q,r}(\a)\ \ch^{SO(2n+1)}_{(2^p,1^q,0^r)}(\x,\ov\x,1)\,,
\end{equation}
are specified for various $a_0$ in Table~\ref{Tab-oon-psi-A11}.
Unlike the situation for $Sp(2n)$, changing the sign of $a_1$ does not result merely in an additional factor 
of $(-1)^q$. One has to simultaneously modify $\psi_{q,r}$ more substantially, as is shown in 
the examples of Table~\ref{Tab-oon-psi-A-11}. 

For $O(2n)$ Corollary~\ref{Cor-eon-AB-AC} allows one to evaluate the expansion
\begin{equation}\label{eqn-eon-psi-A01}
\prod_{i=1}^n (a_0+x_i^2+\ov{x}^2) = \sum_{p=0}^n \sum_{q=0}^{n-p} \ \delta_{r,n-p-q}\ a_2^p\,\psi_{q,r}(\a)\ \ch^{O(2n)}_{(2^p,1^q,0^r)}(\x,\ov\x)\,.
\end{equation}
in the case $m=2$ with $\a=(a_0,0,1)$. As for $Sp(2n)$ the coefficients may only be non-zero if $q=2k$ is even. 
Then setting $q=2k$ yields for example the results for various $\a=(a_0,0,1)$ displayed in Table~\ref{Tab-eon-psi-A01}

Finally, from Theorem~\ref{The-eon} we may evaluate the coefficients $\psi_{q,r}(\a)$ in the expansion
\begin{equation}\label{eqn-eon-psi-A11}
\prod_{i=1}^n (a_0+x_i+\ov{x}_i+x_i^2+\ov{x}^2) = \sum_{p=0}^n \sum_{q=0}^{n-p} \ \delta_{r,n-p-q}\ a_2^p\, \psi_{q,r}(\a)\ \ch^{O(2n)}_{(2^p,1^q,0^r)}(\x,\ov\x)\,.
\end{equation}
Some of their values for various $\a=(a_0,1,1)$ are displayed in Table~\ref{Tab-eon-psi-A11}.


\begin{table}[ht]
\begin{center}
\begin{tabular}{|c||l|l|}
\hline
$\a$&$\psi_{2k,r}(\a)~~\mbox{with}~~\psi_{2k+1,r}(\a)=0$&\cr
\hline
\hline
$(0,0,1)$&$\begin{array}{l}
			 \psi_{2k+8,r}=\psi_{2k,r+4}=\psi_{2k,r}~~\mbox{and}\cr
			 \psi_{2k,r}=\mbox{\begin{footnotesize}$\begin{array}{|c||c|c|c|c|}
                      \hline
                      k\backslash r&0&1&2&3\cr
                      \hline\hline                      
											0&1&-1&0&0\cr
                       \hline
											1&0&1&-1&0\cr
                       \hline
                      2&-1&1&0&0\cr
                       \hline
                      3&0&-1&1&0\cr
                       \hline
                      \end{array}$\end{footnotesize}}\cr\cr
\end{array}$
&\begin{footnotesize}$\begin{array}{c}\mbox{\cite{LO1} $\psi_4(k,r)$}\cr\mbox{in~Table~5.3}\cr\end{array}$\end{footnotesize}\cr
\hline
$(1,0,1)$&$\begin{array}{l}
       \psi_{2k,r}~\mbox{where for all $k,r\geq0$}\cr
       \psi_{2k+6,r}=-\psi_{2k,r+3}=\psi_{2k,r}~~\mbox{and}\cr
			 \psi_{2k,r}=\mbox{\begin{footnotesize}$\begin{array}{|c||c|c|c|}
                      \hline
                      k\backslash r&0&1&2\cr
                      \hline\hline                      0&1&0&0\cr
                       \hline
											1&-1&1&0\cr
                       \hline
                      2&0&-1&0\cr
											\hline
											\end{array}$\end{footnotesize}}\cr\cr
\end{array}$
&\begin{footnotesize}$\begin{array}{c}\mbox{\cite{LO1}~$\psi_3(k,r)$}\cr\mbox{in~Table~5.3}\cr\end{array}$\end{footnotesize}\cr
\hline
$(\ov1,0,1)$&$\begin{array}{l}
       \psi_{2k,r}~\mbox{where for all $k,r\geq0$}\cr
       \psi_{2k,r+6}=-\psi_{2k+6,r}=\psi_{2k,r}~~\mbox{and}\cr
			 \psi_{2k,r}=\mbox{\begin{footnotesize}$\begin{array}{|c||c|c|c|c|c|c|}
                      \hline
                      k\backslash r&0&1&2&3&4&5\cr
                      \hline\hline                      0&1&-2&2&-1&0&0\cr
                       \hline
											1&1&-1&0&1&-1&0\cr
                       \hline
                      2&0&1&-2&2&-1&0\cr
											 \hline
											\end{array}$\end{footnotesize}}\cr\cr
\end{array}$
&\begin{footnotesize}$\begin{array}{c}\mbox{\cite{LO1} $\psi_6(k,r)$}\cr\mbox{in~Table~5.3}\cr\end{array}$\end{footnotesize}\cr
\hline
$(2,0,1)$&$\begin{array}{ll}
           (-1)^k\,(2k+r+2)/2&\mbox{if $r=0 \mod 2$}\cr
			     (-1)^k\, (r+1)/2&\mbox{if $r=1 \mod 2$}\cr
						\end{array}$
&\begin{footnotesize}$\begin{array}{c}\mbox{\cite{LO1} $\psi_2(k,r)$}\cr\mbox{(5.9)}\cr\end{array}$\end{footnotesize}\cr
\hline
$(\ov2,0,1)$&$(-1)^r\, (r+1)(2k+r+2)/2$   
&\begin{footnotesize}$\begin{array}{c}\mbox{\cite{LO1} $\psi_1(k,r)$}\cr\mbox{(5.8)}\cr\end{array}$\end{footnotesize}\cr
\hline
$(3,0,1)$&$(-1)^k\, F_{r+1}\, F_{2k+r+2}$&\cr  
\hline
$(\ov3,0,1)$&$(-1)^r\, \left(\,F_{2k+2r+3)}-F_{2k+1}\,\right)$&\cr  
\hline
\end{tabular}
\end{center}
\medskip
\caption{The $Sp(2n)$ coefficients $\psi_{2k,r}(\a)$ for various $\a=(a_0,0,1)$} 
\label{Tab-spn-psi-A01}																																								
\end{table}

\begin{table}[ht]
\begin{center}
\begin{tabular}{|c||l|}
\hline
$\a$&$\psi_{q,r}(\a)$\cr
\hline\hline$(0,1,1)
$&$\begin{array}{l}\psi_{q,r}~\mbox{for all $q,r\geq0$ with $\psi_{q,r}=$}\cr
\mbox{\begin{footnotesize}$
\begin{array}{|c||c|c|c|}
\hline
q\backslash r&0\mod3&1\mod3&2\mod3\cr
\hline\hline
0\mod3&(q+2r+3)/3&-(q+r+2)/3&-(r+1)/3\cr
\hline
1\mod3&-(q+r+2)/3&-(q+2r+3)/3&(r+1)/3\cr
\hline
2\mod3&(q+1)/3&-(q+1)/3&0\cr
\hline
\end{array}$\end{footnotesize}}\cr\cr
\end{array}$\cr
\hline
$(1,1,1)$&$\begin{array}{l}
\psi_{q,r}~\mbox{for all $q,r\geq0$}\cr
\psi_{q,r+5}=-\psi_{q+5,r}=\psi_{q,r}~~\mbox{and}\cr
\psi_{q,r}=\mbox{\begin{footnotesize}$
\begin{array}{|c||c|c|c|c|c|}
\hline
q\backslash r&0&1&2&3&4\cr
\hline\hline0&1&0&-1&0&0\cr
\hline
1&1&-1&0&0&0\cr
\hline
2&0&0&1&-1&0\cr
\hline
3&0&1&0&-1&0\cr
\hline
4&0&0&0&0&0\cr
\hline
\end{array}$\end{footnotesize}}\cr\cr
\end{array}$\cr
\hline
$(2,1,1)$&$\begin{array}{l}
\psi_{q,r}~\mbox{for all $p,q,r\geq0$}\cr
\psi_{q+12,r}=-\psi_{q,r+6}=\psi_{q,r}~~\mbox{and}\cr
\psi_{q,r}=\mbox{\begin{footnotesize}$
\begin{array}{|c||c|c|c|c|c|c|}
\hline
q\backslash r&0&1&2&3&4&5\cr
\hline\hline0&1&1&1&1&0&0\cr
\hline
1&1&0&1&0&0&0\cr
\hline
2&-1&-1&0&-1&1&0\cr
\hline
3&-2&0&-1&0&1&0\cr
\hline
4&0&1&-1&1&-1&0\cr
\hline
5&2&0&0&0&-2&0\cr
\hline
6&1&-1&1&-1&0&0\cr
\hline
7&-1&0&1&0&2&0\cr
\hline
8&-1&1&0&1&1&0\cr
\hline
9&0&0&-1&0&-1&0\cr
\hline
10&0&-1&-1&-1&-1&0\cr
\hline
11&0&0&0&0&0&0\cr
\hline
\hline
\end{array}$\end{footnotesize}}\cr\cr
\end{array}$\cr
\hline
\end{tabular}
\end{center}
\medskip
\caption{The $Sp(2n)$ coefficients $\psi_{q,r}(\a)$ in (\ref{eqn-spn-psi-A11}) for $\a=(a_0,1,1)$} 
\label{Tab-spn-psi-A11}																																								
\end{table}

\begin{table}[ht]
\begin{center}
\begin{tabular}{|c||l|}
\hline
$\a$&$c_{q,r}(\a)$\cr
\hline
\hline
$(0,0,1)$&$\begin{array}{l}
       \psi_{q,r}~\mbox{where for all $q,r\geq0$}\cr
			 \psi_{q,r+4}=-\psi_{q+4,r}=\psi_{q,r}~~\mbox{and}\cr
			 \psi_{q,r}=\mbox{\begin{footnotesize}$\begin{array}{|c||c|c|c|c|}
                      \hline
                      q\backslash r&0&1&2&3\cr
                      \hline\hline                      0&1&0&-1&0\cr
                       \hline
											1&-1&1&-1&1\cr
                       \hline
                      2&0&1&0&-1\cr
                       \hline
                      3&0&0&0&0\cr
											\hline
										\end{array}$\end{footnotesize}}\cr\cr
\end{array}$\cr
\hline
$(1,0,1)$&$\begin{array}{l}
       \psi_{q,r}~\mbox{where for all $q,r\geq0$}\cr
			 \psi_{q+6,r}=-\psi_{q,r+3}=\psi_{q,r}~~\mbox{and}\cr
			 \psi_{q,r}=\mbox{\begin{footnotesize}$\begin{array}{|c||c|c|c|}
                      \hline
                      q\backslash r&0&1&2\cr
                      \hline\hline                      
											0&1&1&0\cr
                       \hline
											1&-1&1&-1\cr
                       \hline
                      2&-1&0&1\cr
                       \hline
                      3&1&-1&1\cr
											\hline
											4&0&-1&-1\cr
                       \hline
											5&0&0&0\cr
                       \hline
                      \end{array}$\end{footnotesize}}\cr\cr
\end{array}$\cr
\hline
$(\ov1,0,1)$&$\begin{array}{l}
       \psi_{q,r}~\mbox{where for all $q,r\geq0$}\cr
			 \psi_{q,r+6}=-\psi_{q+6,r}=\psi_{q,r}~~\mbox{and}\cr
			 \psi_{q,r}=\mbox{\begin{footnotesize}$\begin{array}{|c||c|c|c|c|c|c|}
                      \hline
                      q\backslash r&0&1&2&3&4&5\cr
                      \hline\hline                      
											0&1&-1&0&1&-1&0\cr
                       \hline
											1&-1&1&-1&1&-1&1\cr
                       \hline
                      2&1&0&-1&1&0&-1\cr
                       \hline
                      3&-1&1&-1&1&-1&1\cr
											\hline
											4&0&1&-1&0&1&-1\cr
                       \hline
											5&0&0&0&0&0&0\cr
                       \hline
                      \end{array}$\end{footnotesize}}\cr\cr
\end{array}$\cr
\hline
$(2,0,1)$&$\begin{array}{cl} (-1)^{q/2}\,(q+2r+2)/2&\mbox{if $q=0 \mod 2$}\cr (-1)^{(q+2r+1)/2}\,(q+1)/2&\mbox{if $q=1 \mod 2$}\cr\end{array}$\cr
\hline
$(\ov2,0,1)$&$\begin{array}{cl} (-1)^{r}\,(q+2r+2)/2&\mbox{if $q=0 \mod 2$}\cr (-1)^{r+1}\,(q+1)/2&\mbox{if $q=1\mod 2$}\cr\end{array}$\cr
\hline
$(3,0,1)$&$\begin{array}{cl} (-1)^{q/2}\,F_{q+2r+2}&\mbox{if $q=0 \mod 2$}\cr (-1)^{(q+2r+1)/2}\,F_{q+1}&\mbox{if $q=1 \mod 2$}\cr\end{array}$\cr
\hline
$(\ov3,0,1)$&$\begin{array}{cl} (-1)^{r}\,F_{q+2r+2}&\mbox{if $q=0 \mod 2$}\cr (-1)^{r+1}\,F_{q+1}&\mbox{if $q=1 \mod 2$}\cr\end{array}$\cr
\hline
\end{tabular}
\end{center}
\medskip
\caption{The $SO(2n+1)$ coefficients $\psi_{q,r}(\a)$ in (\ref{eqn-oon-psi-A01}) for various $\a=(a_0,0,1)$} 
\label{Tab-oon-psi-A01}																																								
\end{table}	

\begin{table}[ht]
\begin{center}
\begin{tabular}{|c||l|}
\hline
$\a$&$\psi_{q,r}(\a)$\cr
\hline\hline$(0,1,1)$&$\begin{array}{l}
\psi_{q,r}~\mbox{for all $q,r\geq0$ with}\cr 
\mbox{\begin{footnotesize}$
\psi_{q,r}=
\begin{array}{|c||c|c|c|}
\hline
q\backslash r&0&1&2\cr
\hline\hline0&1&-1&0\cr
\hline
1&0&-1&1\cr
\hline
2&0&0&0\cr
\hline
\end{array}$\end{footnotesize}}\cr\cr
\end{array}$\cr
\hline
$(1,1,1)$&$\begin{array}{l}
\psi_{q,r}~\mbox{for all $q,r\geq0$}\cr
-\psi_{q+5,r}=\psi_{q,r+5}=\psi_{q,r}~~\mbox{and}\cr
\psi_{q,r}=\mbox{\begin{footnotesize}$
\begin{array}{|c||c|c|c|c|c|}
\hline
q\backslash r&0&1&2&3&4\cr
\hline\hline0&1&0&0&-1&0\cr
\hline
1&0&0&0&-1&1\cr
\hline
2&-1&1&0&0&0\cr
\hline
3&0&1&0&0&-1\cr
\hline
4&0&0&0&0&0\cr
\hline
\end{array}$\end{footnotesize}}\cr\cr
\end{array}$\cr
\hline
$(2,1,1)$&$\begin{array}{l}
\psi_{q,r}~\mbox{for all $q,r\geq0$}\cr
\psi_{q+12,r}=-\psi_{q,r+6}=\psi_{q,r}~~\mbox{and}\cr
\psi_{q,r}=\mbox{\begin{footnotesize}$
\begin{array}{|c||c|c|c|c|c|c|}
\hline
q\backslash r&0&1&2&3&4&5\cr
\hline\hline0&1&1&2&1&1&0\cr
\hline
1&0&1&1&0&1&-1\cr
\hline
2&-2&0&-2&0&0&0\cr
\hline
3&-1&-1&-2&1&-1&2\cr
\hline
4&2&-1&1&0&-1&1\cr
\hline
5&2&0&2&-2&0&-2\cr
\hline
6&-1&1&0&-1&1&-2\cr
\hline
7&-2&1&-1&2&1&1\cr
\hline
8&0&0&0&2&0&2\cr
\hline
9&1&-1&0&-1&-1&0\cr
\hline
10&0&-1&-1&-2&-1&-1\cr
\hline
11&0&0&0&0&0&0\cr
\hline
\end{array}$\end{footnotesize}}\cr\cr
\end{array}$\cr
\hline
\end{tabular}
\end{center}
\medskip
\caption{The $SO(2n+1)$ coefficients $\psi_{q,r}(\a)$ in (\ref{eqn-oon-psi-A11}) for various $\a=(a_0,1,1)$} 
\label{Tab-oon-psi-A11}																																								
\end{table}

\begin{table}[ht]
\begin{center}
\begin{tabular}{|c||l|}
\hline
$\a$&$\psi_{q,r}(\a)$\cr
\hline\hline
$(0,\ov1,1)$&$\begin{array}{l}
\psi_{q,r}~\mbox{for all $q,r\geq0$}\cr
\mbox{\begin{footnotesize}$
\begin{array}{|c||c|c|c|}
\hline
q\backslash r&0\mod3&1\mod3&2\mod3\cr
\hline\hline
0\mod3&(2q\!+\!2r\!+\!3)/3&(2r\!+\!1)/3&(2q\!+\!4r\!+\!4)/3\cr
\hline
1\mod3&-(2q\!+\!4r\!+\!3)/3&-(2r\!+\!1)/3&-(2q\!+\!2r\!+\!3)/3\cr
\hline
2\mod3&(2q\!+\!2)/3&0&(2q\!+\!2)/3\cr
\hline
\end{array}$\end{footnotesize}}\cr\cr
\end{array}$\cr
\hline
$(1,\ov1,1)$&$\begin{array}{l}
\psi_{q,r}~\mbox{for all $q,r\geq0$}\cr
\psi_{q+5,r}=\psi_{q,r+5}=\psi_{q,r}~~\mbox{and}\cr
\psi_{q,r}=\mbox{\begin{footnotesize}$
\begin{array}{|c||c|c|c|c|c|}
\hline
q\backslash r&0&1&2&3&4\cr
\hline\hline0&1&2&-2&-1&0\cr
\hline
1&-2&0&2&-1&1\cr
\hline
2&1&-1&2&0&-2\cr
\hline
3&0&-1&-2&2&1\cr
\hline
4&0&0&0&0&0\cr
\hline
\end{array}$\end{footnotesize}}\cr\cr
\end{array}$\cr
\hline
$(2,\ov1,1)$&$\begin{array}{l}
\psi_{q,r}~\mbox{for all $q,r\geq0$}\cr
\psi_{q+12,r}=-\psi_{q,r+6}=\psi_{q,r}~~\mbox{and}\cr
\psi_{q,r}=\mbox{\begin{footnotesize}$
\begin{array}{|c||c|c|c|c|c|c|}
\hline
q\backslash r&0&1&2&3&4&5\cr
\hline\hline0&1&3&2&3&1&0\cr
\hline
1&-2&-1&-1&-2&1&-1\cr
\hline
2&0&-4&0&-2&0&2\cr
\hline
3&3&3&0&3&-3&0\cr
\hline
4&-2&3&-1&0&1&-3\cr
\hline
5&-2&-4&2&-2&4&2\cr
\hline
6&3&-1&0&1&-3&2\cr
\hline
7&0&3&-3&0&-3&-3\cr
\hline
8&-2&0&2&0&4&0\cr
\hline
9&1&-1&2&1&1&2\cr
\hline
10&0&-3&-2&-2&-3&-1\cr
\hline
11&0&0&0&0&0&0\cr
\hline
\end{array}$\end{footnotesize}}\cr\cr
\end{array}$\cr
\hline
\end{tabular}
\end{center}
\medskip
\caption{The $SO(2n+1)$ coefficients $\psi_{q,r}(\a)$ in (\ref{eqn-oon-psi}) for $\a=(a_0,\ov1,1)$} 
\label{Tab-oon-psi-A-11}																																								
\end{table}


\begin{table}[ht]
\begin{center}
\begin{tabular}{|c||l|l|}
\hline
$\a$&$\psi_{2k,0}(\a)$&$\psi_{2k,r}(\a)$~~\mbox{with $r\geq1$}\cr
\hline
$(0,0,1)$&
$\begin{array}{l}
    \psi_{k,0}~~\mbox{for all $k\geq0$ with}\cr
    \psi_{k+4,0}=\psi_{k,0}~~\mbox{and}\cr
		\psi_{k,0}=\mbox{\begin{footnotesize}$\begin{array}{|c||r|}
                      \hline
                      k\backslash r&0\cr
                      \hline\hline                      0&1\cr
                       \hline
											1&-1\cr
                       \hline
                      2&-1\cr
											\hline
											3&1\cr
											\hline
											\end{array}$\end{footnotesize}}\cr\cr
\end{array}$
&$\begin{array}{l}
     \psi_{k,r}~~\mbox{for all $k\geq0$ and $r\geq1$}\cr
     \psi_{k,r+4}=-\psi_{k+2,r}=\psi_{k,r}~~\mbox{and}\cr
		 \psi_{k,r}=\mbox{\begin{footnotesize}$\begin{array}{|c||r|r|r|r|}
                      \hline
                      k\backslash r&1&2&3&4\cr
                      \hline\hline                      0&0&0&-2&2\cr
                       \hline
											1&2&0&0&-2\cr
                       \hline
											\end{array}$\end{footnotesize}}\cr\cr
\end{array}$
\cr
\hline
$(1,0,1)$
&$\begin{array}{l}
\psi_{k,0}~~\mbox{for all $k\geq1$ and $r\geq0$}\cr
\psi_{k+3,0}(a)=\psi_{k,0}(a)\cr
			 \psi_{k,r}(a)=\mbox{\begin{footnotesize}$\begin{array}{|c||r|}
                      \hline
                      k\backslash r&0\cr
                      \hline\hline                      
											0&1\cr
                       \hline
											1&-2\cr
                       \hline
                      2&1\cr
											 \hline
											\end{array}$\end{footnotesize}}\cr\cr
\end{array}$
&$\begin{array}{l}
        \psi_{k,r}~~\mbox{for all $k\geq0$ and $r\geq1$ with}\cr
        \psi_{k+3,r}=-\psi_{k,r+3}=\psi_{k,r}~~\mbox{and}\cr
		 	  \psi_{k,r}=\mbox{\begin{footnotesize}$\begin{array}{|c||r|r|r|}
                      \hline
                      k\backslash r&1&2&3\cr
                      \hline\hline                      
											0&1&2&-2\cr
                       \hline
											1&1&-1&4\cr
                       \hline
                      2&-2&-1&-2\cr
											 \hline
											\end{array}$\end{footnotesize}}\cr\cr
\end{array}$\cr
\hline
$(\ov1,0,1)$
&$\begin{array}{l}
      \psi_{k,0}~~\mbox{for all $k\geq0$ with}\cr
      \psi_{k+3,0}=-\psi_{k,0}~~\mbox{and}\cr
			\psi_{k,0}=\mbox{\begin{footnotesize}$\begin{array}{|c||r|}
                      \hline
                      k\backslash r&0\cr
                      \hline\hline                      
											0&1\cr
                       \hline
											1&0\cr
                       \hline
                      2&-1\cr
											 \hline
											\end{array}$\end{footnotesize}}\cr\cr
\end{array}$
&$\begin{array}{l}
      \psi_{k,r}~~\mbox{for all $k\geq0$ and $r\geq1$}\cr
      \psi_{k,r+6}=-\psi_{k+3,r}=\psi_{k,r}~~\mbox{and}\cr
			\psi_{k,r}=\mbox{\begin{footnotesize}$\begin{array}{|c||r|r|r|r|r|r|}
                      \hline
                      k\backslash r&1&2&3&4&5&6\cr
                      \hline\hline                      0&-1&0&0&1&-2&2\cr
                       \hline
											1&1&-1&0&1&-1&0\cr
                       \hline
                      2&2&-1&0&0&1&-2\cr
											 \hline
											\end{array}$\end{footnotesize}}\cr\cr
\end{array}$
\cr
\hline
$(2,0,1)$&$(-1)^k\,(2k+1)$&$
\begin{array}{ll}
       (-1)^k\,(4k+2r+2)&r=0\mod 2\cr
			 (-1)^k\,2r&r=1\mod2\cr
\end{array}$
\cr
\hline
$(\ov2,0,1)$&$1$&$(-1)^r\,2$\cr
\hline
$(3,0,1)$&$(-1)^k\,(F_{2k+2}+F_{2k})$&$\begin{array}{ll}(-1)^k(F_{2k+2r+2}+F_{2k+2r})\cr~~~~+(-1)^{k+r}\,(F_{2k+2}+F_{2k})\end{array}$ \cr
\hline
$(\ov3,0,1)$&$F_{2k+1}$&$(-1)^r\,(F_{2k+2r+1}+F_{2k+1})$\cr
\hline
\end{tabular}
\end{center}
\medskip
\caption{The $O(2n)$ coefficients $\psi_{2k,r}(\a)$ in (\ref{eqn-eon-psi-AC}) for various $\a=(a_0,0,1)$} 
\label{Tab-eon-psi-A01}																																								
\end{table}	


\begin{table}[ht]
\begin{center}
\begin{tabular}{|c||l|l|}
\hline
$\a$&$\psi_{q,0}(\a)$&$\psi_{q,r}(\a)$~~\mbox{with $r\geq1$}\cr
\hline
$(0,1,1)$&
$\begin{array}{l}
    \psi_{q,0}~~\mbox{for all $q\geq0$ with}\cr
    \psi_{q+3,0}=\psi_{q,0}~~\mbox{and}\cr
		\psi_{q,0}=\mbox{\begin{footnotesize}$\begin{array}{|c||r|}
                      \hline
                      q\backslash r&0\cr
                      \hline\hline                      
											0&1\cr
                       \hline
											1&1\cr
                       \hline
                      2&0\cr
											\hline
											\end{array}$\end{footnotesize}}\cr\cr
\end{array}$
&$\begin{array}{l}
     \psi_{q,r}~~\mbox{for all $q\geq0$ and $r\geq1$}\cr
     \psi_{q+3,r}=\psi_{q,r+3}=\psi_{q,r}~~\mbox{and}\cr
		 \psi_{q,r}=\mbox{\begin{footnotesize}$\begin{array}{|c||r|r|r|r|}
                      \hline
                      q\backslash r&1&2&3\cr
                      \hline\hline                      
											0&0&-2&2\cr
                       \hline
											1&-2&0&2\cr
                       \hline
                      2&0&0&0\cr
											\hline
											\end{array}$\end{footnotesize}}\cr\cr
\end{array}$
\cr
\hline
$(1,1,1)$
&$\begin{array}{l}
\psi_{q,0}~~\mbox{for all $q\geq1$ and $r\geq0$}\cr
\psi_{q+5,0}(a)=\psi_{q,0}(a)\cr
			 \psi_{q,0}(a)=\mbox{\begin{footnotesize}$\begin{array}{|c||r|}
                      \hline
                      q\backslash r&0\cr
                      \hline\hline                      
											0&1\cr
                       \hline
											1&1\cr
                       \hline
                      2&-1\cr
											 \hline
											3&-1\cr
											\hline
											4&0\cr
											\hline
											\end{array}$\end{footnotesize}}\cr\cr
\end{array}$
&$\begin{array}{l}
        \psi_{q,r}~~\mbox{for all $q\geq0$ and $r\geq1$ with}\cr
        -\psi_{q+5,r}=\psi_{q,r+5}=\psi_{q,r}~~\mbox{and}\cr
		 	  \psi_{q,r}=\mbox{\begin{footnotesize}$\begin{array}{|c||r|r|r|r|r|}
                      \hline
                      q\backslash r&1&2&3&4&5\cr
                      \hline\hline                      0&1&0&-1&-2&2\cr
                       \hline
											1&-1&1&-2&0&2\cr
                       \hline
                      2&0&2&-1&1&-2\cr
											 \hline
											3&2&1&0&-1&-2\cr
                       \hline
                      4&0&0&0&0&0\cr
											 \hline
											\end{array}$\end{footnotesize}}\cr\cr
\end{array}$\cr
\hline
$(2,1,1)$
&$\begin{array}{l}
      \psi_{q,0}~~\mbox{for all $q\geq0$ with}\cr
      \psi_{q+12,0}=\psi_{q,0}~~\mbox{and}\cr
			\psi_{q,0}=\mbox{\begin{footnotesize}$\begin{array}{|c||r|}
                      \hline
                      q\backslash r&0\cr
                      \hline\hline                      0&1\cr
                       \hline
											1&1\cr
                       \hline
                      2&-2\cr
											 \hline
											3&-3\cr
                       \hline
											4&1\cr
                      \hline
                      5&4\cr
											\hline
											6&1\cr
                       \hline
											7&-3\cr
                       \hline
                      8&-2\cr
											 \hline
											9&1\cr
                       \hline
											10&1\cr
                      \hline
                      11&0\cr
											\hline	
											\end{array}$\end{footnotesize}}\cr\cr
\end{array}$
&$\begin{array}{l}
      \psi_{q,r}~~\mbox{for all $q\geq0$ and $r\geq1$}\cr
      \psi_{q+12,r}=-\psi_{q,r+6}=\psi_{q,r}~~\mbox{and $\psi_{q,r}=$}\cr
			\mbox{\begin{footnotesize}$\begin{array}{|c||r|r|r|r|r|r|}
                      \hline
                      q\backslash r&1&2&3&4&5&6\cr
                      \hline\hline                      
											0&2&4&4&2&2&-2\cr
                       \hline
											1&0&4&0&2&0&-2\cr
                       \hline
                      2&-2&-2&-4&2&-2&4\cr
											 \hline
											3&0&-6&0&0&0&6\cr
                       \hline
											4&2&-2&4&-4&2&-2\cr
                       \hline
                      5&0&4&0&-4&0&-8\cr
											\hline
											6&-2&4&-4&2&-2&-2\cr
                       \hline
											7&0&0&0&6&0&6\cr
                       \hline
                      8&2&-2&4&2&2&4\cr
											 \hline
											9&0&-2&0&-4&0&-2\cr
                       \hline
											10&-2&-2&-4&-4&-2&-2\cr
                       \hline
                      11&0&0&0&0&0&0\cr
											\hline
											\end{array}$\end{footnotesize}}\cr\cr
\end{array}$
\cr
\hline
\end{tabular}
\end{center}
\medskip
\caption{The $SO(2n)$ coefficients $\psi_{q,r}(\a)$ in (\ref{eqn-eon-psi-A11}) for various $\a=(a_0,1,1)$} 
\label{Tab-eon-psi-A11}																																								
\end{table}	

\section{Conclusion}~\label{sec-conclusion}

The derivation of some known identities in the form of expansions of products as sums of Schur functions
of $n$ indeterminates $\x=(x_1,x_2,\ldots,x_n)$ has been recast in a Lie group theoretic framework. This has enabled many 
analogous identities to be established for each of the classical Lie groups of rank $n$ for all values of $n$. 
The approach adopted, based on Proposition~\ref{Pro-WG-invariance} and the Weyl group action on arbitrary weight vectors 
$\kappa$ specified in Table~\ref{Tab-wkappa}, has been shown to lead to successively more complicated recurrence relations 
with respect to the highest weight parameters of the characters of irreducible representations of the Lie groups 
$G=GL(n)$, $Sp(2n)$, $SO(2n+1)$ and $SO(2n)$. 
This increasing complexity is due to the different nature of the $n$th generator of the Weyl group, $W_G$,
of the particular Lie group in question, with the number of distinct elementary subsequence transformations of $\kappa$
being unbounded but countably infinite, for all but $GL(n)$, as evidenced by the necessity of the parameter $t$ taking on
all positive integer values in Tables~\ref{Tab-seq-spn}-\ref{Tab-seq-eon}. Hitherto, this phenomenon appears to have
been unrecognised. Nonetheless, the process has been shown to be tractable for the expansion of Weyl group symmetric products
specified by sequences of parameters $\a=(a_0,a_1)$ and $\a=(a_0,a_1,a_2)$ for $Sp(2n)$, $SO(2n+1)$ and $O(2n)$, 
just as it was parameters $\a=(a_0,a_1,a_2,a_3)$ in the case of $GL(n)$.
The same approach has also been applied to expansions in terms of spin characters of $SO(2n+1)$ and $O(2n)$
of certain Weyl group symmetric products parametrised by $\a=(a_\1,a_3,a_\5)$ with either $a_\5=0$ or $\a_3=0$.

From the recurrence relations obtained in this way, results for various 
specific values of $\a$ have been exhibited, concentrating on those examples leading to coefficients
that are periodic in the parameters $p,q,r$ appearing in the specification of the highest weights.
An alternative approach based on the existence of dual group pair identities (\ref{eqn-dual-ids}) has been described. 
In this approach the periodic nature of the resulting coefficients comes as no surprise since the underlying
parameters are all of the form $e^{i\pi/k}$ or $e^{i2\pi/k}$ for integer $k$. Although not leading to recurrence relations
in any natural way, the dual group pair approach has the merit of providing a check on the completeness
of the recurrence relations and of showing that expansion coefficients in the $SO(2n+1)$ 
and $O(2n)$ spin character identities may to be obtained rather more simply from those of $Sp(2n)$ and 
non-spin $SO(2n+1)$ character identities, respectively.

\appendix

\section{Appendix}\label{Appendix}

Here we present the derivation of the recurrence relations for the coefficients
appearing in the $SO(2n+1)$ and $O(2n)$ spin character identities based on the $m=2$ products
parametrised by $\a=(a_\1,0,a_\5)$ that were deferred from Section~\ref{sec-dual-pair}.

\begin{Theorem}\label{The-oon-spin-gf-AC} In the case $\a=(a_\1,0,a_\5)$
\begin{equation}\label{eqn-oon-spin-gf-AC}
    \prod_{i=1}^n \left(a_{\1}(x_i^{\1}\!+\!x_i^{-\1})+a_{\5}(x_i^{\5}\!+\!x_i^{-\5})\right) 
		= \sum_{p=0}^n\,\sum_{q=0}^{n-p}\, a_{\5}^q\, \psi_{q,r}(\a)\, \ch^{SO(2n+1)}_{(a_{\5}^p,a_{\3}^q,a_{\1}^r)}(\x,\ov\x,1)\,,
\end{equation}
where $n=p+q+r$ and 
\begin{equation}\label{eqn-oon-spin-psi-AC}
\begin{array}{l}
\psi_{q,r}(\a) =\ds \chi_{q,r} - a_{\5}^2\,\chi_{q,r-2} + \sum_{s=2}^r \eta_s\,a_{\1}\,a_{\5}^{s-1}\,\chi_{q,r-s} 
 + \zeta_r\,a_{\5}^{r+1}\,\chi_{q-1,0} + \omega_r\,a_{\5}^{r+2}\,\chi_{q-2,0}\,,    \cr
\mbox{with}\cr
\eta_s\!=\!\begin{cases} 1&\mbox{if $s\!=\!0,2\mod 5$};\cr \!-1&\mbox{if $s\!=\!3,4\mod 5$};\cr 0&\mbox{otherwise},\cr \end{cases}
~\zeta_r\!=\!\begin{cases} 1&\mbox{if $r\!=\!1\mod 5$};\cr \!-1&\mbox{if $r\!=\!0\mod 5$};\cr 0&\mbox{otherwise},\cr \end{cases}
~\omega_r\!=\!\begin{cases} 1&\mbox{if $r\!=\!2\mod 5$};\cr \!-1&\mbox{if $r\!=\!3\mod 5$};\cr 0&\mbox{otherwise},\cr \end{cases}\cr   
\end{array}
\end{equation}
and
\begin{equation}\label{eqn-oon-spin-chi-AC}
\begin{array}{l}
\ds \chi_{q,r}(\a)=Q_qR_r \!+\! \sum_{s=1}^r (\sigma_s\,a_{\1}\,a_{\5}^{s}\, Q_{q-1}R_{r-s} \!+\! \tau_s\,a_{\1}\,a_{\5}^{s+1}\,Q_{q-2}R_{r-s}) \!-\! a_{\5}^4 Q_{q-3}R_{r-1}\cr
\mbox{where}~~
   \sigma_s =\begin{cases} 1&\mbox{if $s\!=\!3\mod 5$};\cr -1&\mbox{if $s\!=\!1\mod 5$};\cr 0&\mbox{otherwise},\cr\end{cases}
	~~~~\mbox{and}~~~
	 \tau_s=\begin{cases} 1&\mbox{if $s\!=\!3\mod 5$};\cr -1&\mbox{if $s\!=\!0\mod 5$};\cr 0&\mbox{otherwise}.\cr\end{cases}\cr
\end{array}
\end{equation}
Here $Q_q=0$ and $R_r=0$ if $q<0$ and $r<0$, respectively, while for $q\geq0$ and $r\geq0$
\begin{equation}\label{eqn-oon-spin-Q-AC}
  Q_0=1~~\mbox{and}~~Q_q=-a_{\1}a_{\5} Q_{q-2}+a_{\1}a_{\5}^2 Q_{q-3}+a_{\5}^5 Q_{q-5}~~\mbox{if $q>0$}\,;
\end{equation}
\begin{equation}
\begin{array}{l}\label{eqn-oon-spin-R-AC}
	R_0=1~~\mbox{and}~~R_r=\ds a_{\1} R_{r-1} - a_{\1} a_{\5}^3 R_{r-4} + a_{\5}^5 R_{r-5}+ \sum_{s=3}^r \rho_s\,a_{\1}^2\,a_{\5}^{s-2}\, R_{r-s}~~\mbox{if $r>0$}\,,\cr
~~~~\mbox{where}~~
   \rho_s=\begin{cases} 1&\mbox{if $s\!=\!2,3\mod 5$};\cr -2&\mbox{if $s\!=\!0\mod 5$};\cr 0&\mbox{otherwise}.\cr\end{cases}\cr
\end{array}
\end{equation}
\end{Theorem}

\noindent{\bf Proof}: 
One can extend (\ref{eqn-oon-spin-wi-AB}) and (\ref{eqn-oon-spin-wn-AB}) as follows with the additional data 
\begin{equation}\label{eqn-oon-spin-wi-AC}
\begin{array}{|l|lllll|}
\hline
(\kappa_i,\kappa_{i+1})&(\ov\5,\ov\3)&(\ov\5,\ov\1)&(\ov\5,\1)&(\ov\5,\3)&(\ov\5,\5)\cr
(\kappa_{i+1}-1,\kappa_{i}+1)&(\ov\5,\ov\3)&(\ov\3,\ov\3)&(\ov\1,\ov\3)&(\1,\ov\3)&(\3,\ov\3)\cr
\hline
(\kappa_i,\kappa_{i+1})&(\ov\3,\5)&(\ov\1,\5)&(\1,\5)&(\3,\5)&\cr
(\kappa_{i+1}-1,\kappa_{i}+1)&(\3,\ov\1)&(\3,\1)&(\3,\3)&(\3,\5)&\cr
\hline
\end{array}
\end{equation}
These lead directly to the transformations of elementary subsequences shown in Table~\ref{Tab-seq-oon-spin-AC}
\begin{table}[ht]
\begin{center}
\begin{tabular}{|l|l|l|}
\hline
$\sigma$&$\tau=w\cdot\sigma$&$\sgn(w)\a(\sigma)$ \cr
\hline
$(\ldots,\5,\ldots)$&$(\ldots,\5,\ldots)$&$+a_{\5}$\cr
$(\ldots,\5,\1,\ldots)$&$(\ldots,\5,\1,\ldots)$&$+a_{\1}a_{\5}$\cr
\hline
$(\ldots,\ov\5,\5,\5,\5,\5,\ldots)$&$(\ldots,\3^5,\ldots)$&$+a_{\5}^5$\cr
$(\ldots,\ov\1,\5,\5,\ldots)$&$(\ldots,\3^3,\ldots)$&$+a_{\1}a_{\5}^2$\cr
$(\ldots,\ov\5,\5,\5,\5.\ldots)$&$(\ldots,\3^3,\1,\ldots)$&$-a_{\5}^4$\cr
$(\ldots,\1,\5,\ldots)$&$(\ldots,\3^2,\ldots)$&$-a_{\1}a_{\5}$\cr
$(\ldots(\ov\5,\5,\5,\ov\5,\5)^t,\1,\5,\ldots)$&$(\ldots,\3^2,\1^{5t},\ldots)$&$-a_{\1}a_{\5}^{5t+1}$\cr
$(\ldots(\ov\5,\5,\5,\ov\5,\5)^t,\ov\5,\5,\5,\ov\1,\5,\ldots)$&$(\ldots,\3^2,\1^{5t+3},\ldots)$&$+a_{\1}a_{\5}^{5t+4}$\cr
$(\ldots(\ov\5,\5,\ov\5,\5,\5,)^t,\ov\1,\5,\ldots)$&$(\ldots,\3,\1^{5t+1},\ldots)$&$-a_{\1}a_{\5}^{5t+1}$\cr
$(\ldots(\ov\5,\5,\ov\5,\5,\5,)^t,\ov\5,\5,\1,\5,\ldots)$&$(\ldots,\3,\1^{5t+3},\ldots)$&$+a_{\1}a_{\5}^{5t+3}$\cr
\hline
$(\ldots,\1,\ldots)$&$(\ldots,\1,\ldots)$&$+a_{\1}$\cr
$(\ldots,\ov\5,\ov\1,\5,\5,\ldots)$&$(\ldots,\1^4,\ldots)$&$-a_{\1}a_{\5}^3$\cr
$(\ldots,\ov\5,\ov\5,\5,\5,\5,\ldots)$&$(\ldots,\1^5,\ldots)$&$+a_{\5}^5$\cr
\hline
$(\ldots,\ov\1,(\ov\5,\5,\ov\5,\5,\5)^t,\ov\1,\5,\ldots)$&$(\ldots,\1^{5t+3},\ldots)$&$+a_{\1}^2a_{\5}^{5t+1}$\cr
$(\ldots,\ov\1,(\ov\5,\5,\ov\5,\5,\5)^t,\ov\5,\5,\1,\5,\ldots)$&$(\ldots,\1^{5t+5},\ldots)$&$-a_{\1}^2a_{\5}^{5t+3}$\cr
$(\ldots,\ov\5,\1,\5,(\ov\5,\5,\ov\5,\5,\5)^t,\ov\1,\5,\ldots)$&$(\ldots,\1^{5t+5},\ldots)$&$-a_{\1}^2a_{\5}^{5t+3}$\cr
$(\ldots,\ov\5,\1,\5,(\ov\5,\5,\ov\5,\5,\5)^t,\ov\5,\5,\1,\5,\ldots)$&$(\ldots,\1^{5t+7},\ldots)$&$+a_{\1}^2a_{\5}^{5t+5}$\cr
\hline
\hline
\end{tabular}
\end{center}
\medskip
\caption{Elementary subsequences $\tau=w\cdot\sigma$ of $\lambda$ along with $\sgn(w)$ and the contribution of 
$a(\sigma)$ to $a(\kappa)$ in the case $SO(2n+1)$, with $t\geq0$.} 
\label{Tab-seq-oon-spin-AC}
\end{table}

Similarly, the additional data from 
\begin{equation}\label{eqn-oon-spin-wn-AC}
\begin{array}{|l|l|}
\hline
(\ldots,\kappa_n)&(\ldots,\ov\5)\cr
(\ldots,-\kappa_n-1)&(\ldots,\3)\cr
\hline
\end{array}
\end{equation}
allows one to derive the terminating subsequence transformations given in Table~\ref{Tab-seq-term-oon-spin-AC}. 
\begin{table}[ht]
\begin{center}
\begin{tabular}{|l|l|l|}
\hline
$\sigma$&$\tau=w\cdot\sigma$&$\sgn(w)\a(\sigma)$ \cr
\hline
\hline
$(\ldots(\ov\5,\5,\5,\ov\5,\5)^t,\ov\5,\5,\5,\ov\5)$&$(\ldots,\3^2,\1^{5t+2})$&$+a_{\5}^{5t+4}$\cr
$(\ldots(\ov\5,\5,\5,\ov\5,\5)^t,\ov\5,\5,\5,\ov\5,\5)$&$(\ldots,\3^2,\1^{5t+3})$&$-a_{\5}^{5t+5}$\cr
$(\ldots(\ov\5,\5,\ov\5,\5,\5)^t,\ov\5)$&$(\ldots,\3,\1^{5t})$&$-a_{\5}^{5t+1}$\cr
$(\ldots(\ov\5,\5,\ov\5,\5,\5)^t,\ov\5,\5)$&$(\ldots,\3,\1^{5t+1})$&$+a_{\5}^{5t+2}$\cr
\hline
$(\ldots,\ov\5,\ov\5,)$&$(\ldots,\1^2)$&$-a_{\5}^2$\cr
\hline
$(\ldots,\ov\1,(\ov\5,\5,\ov\5,\5,\5)^t,\ov\5,)$&$(\ldots,\1^{5t+2})$&$+a_{\1}a_{\5}^{5t+1}$\cr
$(\ldots,\ov\1,(\ov\5,\5,\ov\5,\5,\5)^t,\ov\5,\5)$&$(\ldots,\1^{5t+3})$&$-a_{\1}a_{\5}^{5t+2}$\cr
$(\ldots,\ov\5,\1,\5,(\ov\5,\5,\ov\5,\5,\5,)^t,\ov\5)$&$(\ldots,\1^{5t+4})$&$-a_{\1}a_{\5}^{5t+3}$\cr
$(\ldots,\ov\5,\1,\5,(\ov\5,\5,\ov\5,\5,\5,)^t,\ov\5,\5)$&$(\ldots,\1^{5t+5})$&+$a_{\1}a_{\5}^{5t+4}$\cr
\hline
\hline
\end{tabular}
\end{center}
\medskip
\caption{Elementary subsequences and terminating subsequences $\tau=w\cdot\sigma$ of $\lambda$ along with $\sgn(w)$ and the contribution of 
$a(\sigma)$ to $a(\kappa)$ in the case $SO(2n+1)$, with $t\geq0$.} 
\label{Tab-seq-term-oon-spin-AC}
\end{table}

The repeated use of (\ref{eqn-oon-spin-wi-AB}), (\ref{eqn-oon-spin-wn-AB}), (\ref{eqn-oon-spin-wi-AC}) and (\ref{eqn-oon-spin-wn-AC})
not only allow one to derive all the results in Tables~\ref{Tab-seq-oon-spin-AC} and ~\ref{Tab-seq-term-oon-spin-AC}, 
but also show that all other non-standard subsequences
in $\{\ov\5,\ov\1,\1,\5\}$ lead to contributions that are zero by virtue of changing sign under a Weyl group reflection.
It is perhaps worth pointing out the origin of the $t$-dependent terms in these Tables.  
These subsequences involve either $(\ov\5,\5,\ov\5,\5,\5)^t$ or $(\ov\5,\5,\5,\ov\5,\5)^t$. 
These transform to give $(\3,\1,\1,\1,\ov\1)^t$ and $(\3,\3,\1,\1,\ov\3)^t$, respectively, 
which may be rewritten as $\3,(\1,\1,\1,\ov\1,\3)^{t-1},\1,\1,\1,\ov\1$ and $\3,\3,(\1,\1,\ov\3,\3,\3)^{t-1},\1,\1,\ov\3$.
The subsequences $\ov\1,\3$ and $\ov\3,\3,\3$ transform to give $\1,\1$ and $\1,\1,\1$, respectively, yielding in each case 
a factor of the form $(\1,\1,\1,\1,\1)^{t-1}$, that is $\1^{5(t-1)}$. 
Taking into account further transformation at the head and tail of these subsequences, along
with all the relevant sign changes, one arrives at the $t$-dependent results presented in the two Tables.
It is these tabulations of all possible non-vanishing contributions to the expansion (\ref{eqn-oon-spin-gf-AC})
that complete the proof of the recurrence relations (\ref{eqn-oon-spin-psi-AC})-(\ref{eqn-oon-spin-R-AC}). 
\qed

In the comparable $SO(2n)$ case we find
\begin{Theorem}\label{The-eon-spin-gf-AC} In the case $\a=(a_\1,0,a_\5)$
\begin{equation}\label{eqn-eon-spin-gf-AC}
    \prod_{i=1}^n \left(a_{\1}(x_i^{\1}\!+\!x_i^{-\1})+a_{\5}(x_i^{\5}\!+\!x_i^{-\5})\right) 
		= \sum_{p=0}^n\,\sum_{q=0}^{n-p}\, a_{\5}^q\, \psi_{q,r}(\a)\, \ch^{O(2n)}_{(a_{\5}^p,a_{\3}^q,a_{\1}^r)}(\x,\ov\x)\,,
\end{equation}
where $n=p+q+r$ and
\begin{equation}\label{eqn-eon-spin-psi-AC}
\begin{array}{l}
\psi_{q,r}(\a) =\ds \chi_{q,r} - a_{\5}^3\,\chi_{q,r-3} + \sum_{s=2}^r \eta_s\,a_{\1}\,a_{\5}^{s-1}\,\chi_{q,r-s} 
 + \zeta_r\,a_{\5}^{r+1}\,\chi_{q-1,0} + \omega_r\,a_{\5}^{r+2}\,\chi_{q-2,0}\,, \cr 
\mbox{with}\cr
\eta_s\!=\!\begin{cases} 1&\mbox{if $s\!=\!1,2\mod 5$};\cr \!-1&\mbox{if $s\!=\!0,3\mod 5$};\cr 0&\mbox{otherwise},\cr \end{cases}
~\zeta_r\!=\!\begin{cases} 1&\mbox{if $r\!=\!3\mod 5$};\cr \!-1&\mbox{if $r\!=\!4\mod 5$};\cr 0&\mbox{otherwise},\cr \end{cases}
~\omega_r\!=\!\begin{cases} 1&\mbox{if $r\!=\!1\mod 5$};\cr \!-1&\mbox{if $r\!=\!0\mod 5$};\cr 0&\mbox{otherwise},\cr \end{cases}\cr   
\end{array}
\end{equation}
and with $\chi_{q,r}$ still defined by (\ref{eqn-oon-spin-chi-AC}).
\end{Theorem}

\noindent{\bf Proof}: 
The only difference from the proof of Theorem~\ref{The-oon-spin-gf-AC} is that the terminating subsequence transformations
have to be recalculated using the elementary Weyl group transformations
\begin{equation}\label{eqn-eon-spin-wn-15}
\begin{array}{|l|lllll|}
\hline
(\ldots,\kappa_{n-1},\kappa_n)&(\ldots,\ov\5,\ov\5)&(\ldots,\ov\3,\ov\5)&(\ldots,\ov\1,\ov\5)&(\ldots,\1,\ov\5)&(\ldots,\3,\ov\5)\cr
(\ldots,-\kappa_n\!-\!1,-\kappa_{n-1}\!-\!1)&(\ldots,\3,\3)&(\ldots,\3,\1)&(\ldots,\3,\ov\1)&(\ldots,\3,\ov\3)&(\ldots,\3,\ov\5)\cr
\hline
\end{array}
\end{equation}
These lead directly to the transformations of terminating subsequences shown in Table~\ref{Tab-seq-term-eon-spin-AC}
\begin{table}[ht]
\begin{center}
\begin{tabular}{|l|l|l|}
\hline
$\sigma$&$\tau=w\cdot\sigma$&$\sgn(w)\a(\sigma)$ \cr
\hline
$(\ldots(\ov\5,\5,\5,\ov\5,\5)^t,\ov\5,\ov\5)$&$(\ldots,\3^2,\1^{5t})$&$-a_{\5}^{5t+2}$\cr
$(\ldots(\ov\5,\5,\5,\ov\5,\5)^t,\ov\5,\5,\5,\ov\5,\5)$&$(\ldots,\3^2,\1^{5t+1})$&$+a_{\5}^{5t+5}$\cr
$(\ldots(\ov\5,\5,\ov\5,\5,\5)^t,\ov\5,\5,\ov\5,\ov\5)$&$(\ldots,\3,\1^{5t+3})$&$+a_{\5}^{5t+4}$\cr
$(\ldots(\ov\5,\5,\ov\5,\5,\5)^t,\ov\5,\5,\ov\5,\5,\ov\5)$&$(\ldots,\3,\1^{5t+4})$&$-a_{\5}^{5t+5}$\cr
\hline
$(\ldots,\ov\5,\ov\1)$&$(\ldots,\1^2)$&$+a_{\1}a_{\5}$\cr
$(\ldots,\ov\5,\1,\ov\5)$&$(\ldots,\1^3)$&$-a_{\1}a_{\5}^2$\cr
$(\ldots,\ov\5,\ov\5,\5)$&$(\ldots,\1^3)$&$-a_{\5}^3$\cr
\hline
$(\ldots,\ov\1,(\ov\5,\5,\ov\5,\5,\5)^t,\ov\5,\5,\ov\5,\ov\5)$&$(\ldots,\1^{5t+5})$&$+a_{\1}a_{\5}^{5t+4}$\cr
$(\ldots,\ov\1,(\ov\5,\5,\ov\5,\5,\5)^t,\ov\5,\5,\ov\5,\5,\ov\5)$&$(\ldots,\1^{5t+6})$&$-a_{\1}a_{\5}^{5t+5}$\cr
$(\ldots,\ov\5,\1,\5,(\ov\5,\5,\ov\5,\5,\5,)^t,\ov\5,\5,\ov\5,\ov\5)$&$(\ldots,\1^{5t+7})$&$-a_{\1}a_{\5}^{5t+6}$\cr
$(\ldots,\ov\5,\1,\5,(\ov\5,\5,\ov\5,\5,\5,)^t,\ov\5,\5,\ov\5,\5,\ov\5)$&$(\ldots,\1^{5t+8})$&+$a_{\1}a_{\5}^{5t+7}$\cr
\hline
\hline
\end{tabular}
\end{center}
\medskip
\caption{Terminating subsequences $\tau=w\cdot\sigma$ of $\lambda$ along with $\sgn(w)$ and the contribution of 
$a(\sigma)$ to $a(\kappa)$ in the case $SO(2n)$, with $t\geq0$.} 
\label{Tab-seq-term-eon-spin-AC}
\end{table}
This tabulation leads directly to the required recurrence relations (\ref{eqn-eon-spin-psi-AC}) for $\psi_{q,r}(\a)$ 
in terms of various $\chi_{q,r}$. This completes the proof of (\ref{eqn-eon-spin-gf-AC}) since 
$\chi_{q,r}$, $Q_q$ and $R_r$ are defined by (\ref{eqn-oon-spin-chi-AC}), (\ref{eqn-oon-spin-Q-AC}) and (\ref{eqn-oon-spin-R-AC}). 
\qed

Despite the complexity of the recurrence relations of Theorems~\ref{The-oon-spin-gf-AC} and~\ref{The-eon-spin-gf-AC}, they have been used to 
calculate $\psi_{qr}(\a)$ for both $SO(2n+1)$ and $O(2n)$ in the cases $\a=(a_\1,0,1)$ with $a_\1=0$, $1$ and $-1$. The results are
entirely as predicted by the second two identities of (\ref{eqn-phi-psi-oon-spn}). This is to be seen as evidence not so much for the validity 
of (\ref{eqn-phi-psi-oon-spn}), but as check on the validity of the recurrence relations,


\end{document}